\documentclass[12pt,a4paper]{article}
\usepackage{amsmath,mathtools,pifont,amsfonts,amssymb,latexsym}
\usepackage{theorem}
\usepackage{enumerate,epsfig}
\usepackage{hyperref}
\usepackage{MnSymbol}%for the brokenarrow symbol 

\hypersetup{colorlinks,pdfauthor={Andreas Knauf}}

\numberwithin{equation}{section}
\numberwithin{figure}{section}
\theoremstyle{change}
\newtheorem{theorem}{Theorem}[section]
\newtheorem {lemma}[theorem]{Lemma}
\newtheorem {prop}[theorem]{Proposition}

\newtheorem {corollary}[theorem]{Corollary}
{\theorembodyfont{\normalfont}
  \newtheorem {definition}[theorem]{Definition}
  
  \newtheorem {remark}[theorem]{Remark}
  
  \newtheorem {earlier}[theorem]{Earlier results}
  \newtheorem {caveats}[theorem]{Caveats}
  
  \newtheorem {example}[theorem]{Example}}

\newcommand{\beq}{\begin{equation}}
\newcommand{\eeq}{\end{equation}}
\newcommand{\Leq}[1]{\label{#1}\end{equation}}
\newcommand{\beqn}{\begin{eqnarray}}
\newcommand{\eeqn}{\end{eqnarray}}
\newcommand{\beqno}{\begin{eqnarray*}}
\newcommand{\eeqno}{\end{eqnarray*}}
\renewcommand {\l}{\left}
\newcommand {\ri}{\right}
 
\newcommand {\LA}{\left\langle}
\newcommand {\RA}{\right\rangle}
\newcommand {\pa}{\partial}
\newcommand {\eh}{{\textstyle \frac{1}{2}}}
\newcommand {\ev}{{\textstyle \frac{1}{4}}}
\newcommand {\ea}{{\textstyle \frac{1}{8}}}

\newcommand {\bN}{{\mathbb N}}
\newcommand {\bR}{{\mathbb R}}

\newcommand{\idty}{{\rm 1\mskip-4mu l}} %Identity Operator 
\newcommand{\cC}{{\cal C}} % 
\newcommand{\cD}{{\cal D}} %
\newcommand{\cF}{{\cal F}} % 
\newcommand{\cM}{{\cal M}}
\newcommand{\cO}{{\cal O}} % O symbol
\newcommand{\cP}{{\cal P}}
\newcommand{\ov}{\overline}
\newcommand{\bem}{\l(\! \begin{array}}
\newcommand{\eem}{\end{array}\!\ri)}
\newcommand{\bsm}{\left(\begin{smallmatrix}} % begin small matrix
\newcommand{\esm}{\end{smallmatrix}\right)}  % end   small matrix
\newcommand{\NN}{\nonumber}
\newcommand{\qmbox}[1]{\quad\mbox{#1}\quad}
\renewcommand {\max}{{{\rm max}}}

\begin{document}
\title {Classical $n$-body scattering\\ with long-range potentials}
\date{\today}
\author{Jacques F{\'e}joz\thanks{Ceremade, Universit{\'e} Paris Dauphine --
    PSL \& IMCCE, Observatoire de Paris -- PSL,
    \texttt{jacques.fejoz@dauphine.fr}}
\and
Andreas Knauf\thanks{
Department of Mathematics,
Friedrich-Alexander-University Erlangen-N{\"u}rnberg,
Cauerstr.\ 11, D-91058 Erlangen, 
Germany, \texttt{knauf@math.fau.de}} 
\and
Richard Montgomery\thanks{Mathematics Department,
UC Santa Cruz, 
4111 McHenry, 
Santa Cruz, CA 95064, USA, \texttt{rmont@ucsc.edu}}}
\date{\today}
\maketitle
\begin{abstract}
  We consider the scattering of $n$ classical particles interacting
  via pair potentials, under the assumption that each pair potential
  is ``long-range'', i.e.\ being of order $\cO(r^{-\alpha})$  for some
  $\alpha >0$.  We define and focus on the ``free region'', the set of
  states leading to well-defined and well-separated final states at
  infinity. As a first step, we prove the existence of an explicit,
  global surface of section for the free region. 
  This surface of section  is key to proving  the
  smoothness of the map  sending a point to its final state  and to establishing a 
  forward conjugacy between the $n$-body dynamics and  free dynamics. 
\end{abstract}

\clearpage
\tableofcontents

%----------------------------------------------------
\section{Main results and set-up}
%----------------------------------------------------

%----------------------------------------------------
\subsection{Main results}
%----------------------------------------------------

Consider  $n$ classical particles moving in  $d$-dimensional
Euclidean space under the influence of a potential  which is the
sum of pair potentials.  
If the pair potentials die off
appropriately at infinity then we expect that, within any widely
separated fast-moving configuration of particles, the individual
particles will move almost along straight lines. In this case it makes
sense to talk about ``scattering''.  See for example \cite{DG,Sim,He,Hu1}, and \cite{Hu2}.  
We will prove 
new facts regarding the relation between initial conditions and
scattering data at infinity.  The most surprising of these  is   explicit criteria \eqref{def:finally:free}
which guarantee the  escape to asymptotic freedom.
See Theorem \ref{thm:main1}.   Other facts, summarized by Theorems \ref{thm:main2}, 
\ref{thm:mainMoller}  extend and refine results  
previously only known for short range
potentials  to the case of long-range potentials (see Definition \ref{def:range}).
 
 \subsection{Setup and notation for $n$-body dynamics and potential decay}
\medskip 

\noindent
A configuration $q$ specifies the locations of all $n$
masses, so that $q = (q_1, \ldots, q_n) \in \bR^{dn}$ with
$q_a \in \bR^d$.  Thus our configuration space is \beq M:=\bR^{dn}_q
\qmbox{, or} \widehat M:=\bR^{dn}_q\backslash\Delta,
\Leq{conf:space:hat} depending on whether or not our pair potentials
$V_{i,j} = V_{i,j} (q_i -q_j)$ have singularities at collision
$q_i = q_j$; here \beq \Delta := \{q = (q_1,\ldots,q_n)\in
\bR^{dn}_q\mid q_i=q_j\mbox{ for some }i\neq j\} \Leq{collision:set}
is the collision set, also known as the ``fat diagonal''. $\Delta$
will also play an important role in the velocity space.

Configurations evolve in time according to   Newton's equations
\beq
m_i \ddot q_i = -\nabla_{q_i} V,\qquad i =1, \ldots , n,  \qquad m_i >
0 \text{ the masses},  
\Leq{N}
\noindent which we will  formulate in the usual way in phase space,
using momenta $p_i = m_i \dot q_i$, so that $p = (p_1, \ldots, p_n)
\in \bR^{dn}_p$. Thus our   phase space $P$ is 
\beq
P := T^*M = \bR^{dn}_p \times \bR^{dn}_q, \qmbox{ or} 
\widehat P := T^*\widehat M =  \bR^{dn}_p \times (\bR^{dn}_q \, \backslash\, \Delta) 
\Leq{phase_space}
endowed with its  canonical symplectic form.
Identify $\bR^{dn}$ with $\bR^n \otimes \bR^d$, let 
\[\cM:={\rm diag}(m_1,\ldots,m_n)\otimes \idty_d\]
be the mass matrix, seen as an (invertible symmetric) operator on
$\bR^{dn}_p$.  Newton's equations can be rewritten as Hamilton's
equations
\[\dot p= - \nabla_q V, \quad \dot q = \cM^{-1} p,\]
with Hamiltonian  $H:P\to\bR$  (or $\widehat P \to \bR$),  
\beq 
H(p,q) := K(p) + V(q),
\Leq{Ham}
where the potential energy is assumed to be of the form
\[V(q) := \sum_{1\le i<j\le n} V_{i,j}(q_i-q_j),\]
where the  pair potentials $V_{i,j}$ satisfy  $V_{j,i}=V_{i,j}$
and $V_{i,i}=0$ for all $i,j$, and where $K$ is the usual kinetic energy
$$K(p) := \sum_{i=1}^n \frac{\|p_i\|^2}{2m_i}= \eh \LA
p,p\RA_{\!\cM^{-1}}, \quad \LA p,p'\RA_{\!\cM^{-1}}:= \LA p,\cM^{-1}
p'\RA.$$

\medskip From now on we will use multi-index notation for partial
derivatives.

\begin{definition}
  \label{def:range}%
  A pair potential $V_{i,j}\in C^2(\bR^d \backslash \{0\},\bR)$ is 
  \begin{enumerate}[$\bullet$]
  \item 
    \emph{long range} if for some $\alpha>0$
    \beq
    \partial^\gamma V_{i,j}(q)=\cO \big(\|q\|^{-\alpha -
      |\gamma|}\big) \qquad \big(\|q\| \to \infty,\,
    \gamma\in\bN_0^d,\, |\gamma|\le 2 \big) 
    \Leq{V:decay} 
    (if needed, $V_{i,j}$ will then also be called an $\alpha$-potential),
  \item 
    \emph{short range} if \eqref{V:decay} is valid for some $\alpha>1$,
  \item \emph{finite range} if the $V_{i,j}$ have bounded support.
  \end{enumerate}
  The potential $V(q)=\sum V_{i,j}(q_i-q_j)$ is called \emph{long
    range}, etc., if all its pair potentials $V_{i,j}$ have the
  corresponding property.
\end{definition}

\begin{caveats}\quad\\
 According to this established
  terminology, the following implications hold:
  $$\mbox{finite range } \Longrightarrow \mbox{ short range }
  \Longrightarrow \mbox{ long range}.$$ We apologize for the
  counterintuitive nature of the terminology.  It is standard in
  scattering literature.  Also note that a finite range potential $V$
  typically does not have bounded support within ${\bR}^{dn}$.
  Rather, its support is contained in a neighborhood of the fat
  diagonal $\Delta$.  \hfill $\Diamond$
\end{caveats}

\begin{example}[Celestial mechanics and electrostatics]\quad\\
  \label{ex:grav:coulomb}%
  In celestial mechanics and electrostatics we have
  $V_{i,j}(Q)=\frac{I_{i,j}}{\|Q\|} $ with respectively
  $I_{i,j}=-m_i m_j$ and $I_{i,j}=Z_i Z_j$ for the charges
  $Z_i\in\bR\!\setminus\!  \{0\}$. These potentials are long range,
  lying on the boundary of the space of short range potentials.
  \hfill $\Diamond$
\end{example}

\begin{remark}[Strong forces near collisions]\quad\\
  By definition, so-called strong force potentials satisfy
  \beq
  \partial^\gamma V_{i,j}(q)=\cO \big(\|q\|^{-\alpha - |\gamma|}\big)
  \qquad \big(\|q\| \to 0,\, \gamma\in\bN_0^d,\, |\gamma|\le 2 \big),
  \Leq{}
  for some $\alpha \geq 2$ (cf.~\eqref{V:decay} as
  $q \to \infty$).  Variationally speaking, this condition is most
  important in the opposite ''ultraviolet'' regime of short distances,
  $\|q_i - q_j \| \ll 1$, rather than our current ``infrared regime'' of
  long distances.  Imposing the strong force condition on attractive
  forces guarantees that any collision solution has infinite action
  and so is a simple way to exclude collision solutions as candidate
  minimizers when using the direct method of the calculus of
  variations to achieve various types of solutions (e.g.\ periodic
  ones).\hfill $\Diamond$
\end{remark}

%----------------------------------------------------
\subsection{Asymptotic freedom} 
%----------------------------------------------------

Our first goal is to define the \emph{free region} of phase space,
leading to motions along which mutual distance eventually increase
linearly with time, as in the free flow, where bodies do not interact.

This definition relies on the prior concept of asymptotic velocity. 

\begin{definition}
  \label{def:asymptotic:velocity}%
  The \emph{(forward, resp. backward) asymptotic shape} or
  \emph{velocity} of a state $x \in P$ is the limit in $\bR^{dn}$, if
  it exists, is
  \[v^\pm (x) := \lim_{t \to \pm \infty} \frac{q(t)}{t},\]
  where $x(t)=(q(t),p(t))$ is the integral curve through $x$ at $t=0$.
\end{definition}

We are interested in motions for which $v^+ \notin \Delta$.

\begin{definition}
  \label{def:free}\quad\\%
  The state $x$ is \emph{forward free} if $v^+(x)$ exists and
  $v^+(x) \in \bR^{dn} \, \backslash\, \Delta$. We call
  \[F^+ := \{x \in P\mid v^+(x) \in \bR^{dn} \, \backslash\, \Delta\}\]
  the subset of $P$ of forward free states.  
  Correspondingly, the subset $F^-$ is the set of states $x$ which are
  \emph{backward free}, i.e. $v^-(x) \in \bR^{dn} \, \backslash\, \Delta$.
  
  We will sometimes  refer to trajectories passing through $F^+$
  as \emph{escape orbits}. 
\end{definition} 
 
\begin{remark}[Clusters] \quad\\
 Those  motions $x(t)$ for which  $v^+ (x)$ exists but  for which   $v^+ (x)  \in \Delta$   break up into $k < n$ clusters, each cluster
composed of those particles whose  indices $i$ share a common asymptotic velocity:  $v^+_i = v^+ _j$.
The dynamics within a cluster of size $c$ can be as complicated as that of the general $c$-body problem. The
clusters interact with each other   like a  free 
$k$-body system.  (See \textsc{Marchal-Saari}~\cite{MS}, however not in the sense of asymptotic completeness, see \cite[section 5.10]{DG}.)\,$\Diamond$%
\end{remark} 
 
\begin{example}[Celestial mechanics]\quad\\
  \label{ex:celestial:mechanics}%
  \textsc{Chazy} \cite{Cha} showed
  that collision-free solutions for $n=3$ gravitating bodies fall into one of seven possibilities,
  regarding their final behavior in the future. 
  \begin{itemize}
  \item 
  Bounded, parabolic, parabolic-elliptic and oscillating motions have zero asymptotic velocity.
  \item 
  Hyperbolic-elliptic and hyperbolic-parabolic motions have asymptotic velocity belonging to 
  $\Delta\backslash\{0\}$.
  \item Hyperbolic motions are free. So their asymptotic velocity is in $\bR^{3d}\,\backslash\,\Delta$. 
  \end{itemize}
  So, here hyperbolicity equates to freedom. For more bodies, new types
  of final motions occur, notably the ``non-collision singularities'', see {\sc Gerver} \cite{Ge}
  and {\sc Xia} \cite{Xia}.
  But it remains true that every
  collision-free solution has  asymptotic velocities $v^{\pm}$ in both time directions
  provided we allow  velocities to take values in the one point compactification
  of $\bR^{dn}$. (For example,  for initial conditions $x$ leading to non-collision singularities  we have 
  $\lim_{t \to T^\pm(x)} \frac{\|q(t)\|}{t} =\infty$, where $T^+(x)\in(0,\infty]$ and $T^-(x)\in[-\infty,0)$ are
  the escape times  beyond which the solution fails to exist.)
  \hfill$\Diamond$%
\end{example}

The precise structure of $F^+$ is not obvious. Yet, by flowing $F^+$
along integral curves, we will reach an open subset of $P$, which we can characterize explicitly. Let
\begin{align}
  q_{i,j}& := \|q_i-q_j\|\ ,\quad q_{\min} := \min_{i<j}q_{i,j} 
           \ ,\quad  q_{\max}:=  \operatorname*{\max}_{i<j}
           q_{i,j}\,,\NN\\
  v_{i,j}&:=\|v_i-v_j\|\ ,\quad  
           v_{\min}:=\min_{i<j} v_{i,j}\ ,\quad  
           v_{\max}:= \operatorname*{\max}_{i<j} v_{i,j}
\end{align}
and let  $\alpha$, $\delta$ and $C$ be three positive parameters.

\begin{definition}
  \label{def:finally:free}%
  The \emph{finally free region} (with parameters $\alpha$, $\delta$
  and $C$) is 
  \beqn F^+_{\rm loc} :=
  \Big\{ x=(p,q) \in P\hspace{-2mm} & \Big| \hspace{-2mm}&
  v_{\min}^2> C \frac{q_{\max}}{q_{\min}^{\alpha+1}}\,
  , \label{def:ff}\\
  && \LA v_i-v_j, q_i-q_j\RA > (1-\delta)v_{i,j}q_{i,j}\, ,\NN\\
  &&\textstyle (1+2\delta) \frac{q_{k,l}}{v_{k,l}}>
  \frac{q_{i,j}}{v_{i,j}} ,\quad( i \ne j,k \ne l )\Big\}. \NN
  \eeqn
\end{definition}
Notice that $F^+_{\rm loc}$, like $F^+$, is invariant w.r.t.\ the symplectic lift of the diagonal action
of the Euclidean group on configuration space $\bR^{dn}$.

The following theorem justifies that our definition of $F^+$ matches
our initial goal, and also justifies the notation $F^+_{\rm loc}$.

\begin{theorem} 
  \label{thm:main1}%
  For any long range potential $V$, there exist appropriate parameters
  $\alpha$, $\delta$ and $C$ such that $F^+_{\rm loc}$ is forward invariant and such that a
  state $x \in P$ is in $F^+$ if and only if its forward orbit
  eventually enters $F^+_{\rm loc}$. 
\end{theorem}

The theorem follows from Theorem~\ref{thm:final:free} below. The proof
will actually show that the boundary $\partial F^+_{\rm loc}$ is a ($C^0$)
surface of section of the flow restricted to $F^+$.
Notice that
$F^+ = \bigcup_{t\geq 0} \Phi_{-t}\left(F^+_{\rm loc}\right)$, where
$F^+_{\rm loc}$ is open and $\Phi_{-t}$ is smooth, whence the following.

\begin{corollary}
  $F^+$ is a non-empty open subset of $P$. 
\end{corollary}

\medskip The asymptotic velocity map \ref{def:asymptotic:velocity} 
enjoys regularity on $F^+$.
 
\begin{theorem}
  [Asymptotic velocity map on $F^+$]\quad \label{thm:main2}\\
  Assume that $V$ is a long-range potential whose pair potentials are
  $C^k$, $k \ge 2$. The map
  $v^+ : F^+ \to \bR^{dn} \, \backslash\, \Delta$ is a $C^{k-1}$
  complete set of commuting first integrals. Moreover, for fixed
  $v_* \in \bR^{dn} \, \backslash\, \Delta$ the space of all forward
  orbits $x(t)$ for which $v^+ (x(t)) = v_*$ has the structure of an
  affine space modelled on the $(nd-1)$-dimensional vector
  space~$v_* ^{\perp}$.
\end{theorem} 

The regularity of the $v^+$ follows from item 1 of Theorem
\ref{thm:both:moeller}. That components of $v^+$ Poisson commute is
clear. That $v^+$ is a surjective submersion, and the assertion on the
structure of its fibers follows from item 4 of Theorem
\ref{thm:Dollard:Moeller} below.
\begin{earlier}[Smoothness of scattering data]\quad\\
  Smoothness of the scattering data and in particular of the
  asymptotic velocity map $x \mapsto v^+ (x)$ has been achieved under
  various assumptions:
\begin{enumerate}[$\bullet$]
\item 
In \cite{Gu}, {\sc Gutkin} proved continuity of scattering data for
a class of $n$-particle systems on the line with repulsive interactions.
\item 
Later, {\sc Fusco} and {\sc Oliva} proved in \cite{FO} a result that
implies smoothness of asymptotic momentum and even integrability for
repulsive Coulombic potentials.
\item 
More recently, \textsc{Duignan} et al. \cite{Duignan} prove that the
map $x \to v^+(x)$ is analytic on $F^+_{\rm loc}$ for the Newtonian
potential. \hfill $\Diamond$
\end{enumerate}
\end{earlier}

%----------------------------------------------------
\subsection{Comparison with free flows}
\label{subsec:flows}
%----------------------------------------------------

In order to study the asymptotic behavior of the dynamics on $F^+$,
one strategy would be to compactify the phase space, as
in~\cite{Duignan} for the $N$-body problem. Such a
compactification is hard to define in full generality. Another
strategy, chosen here, is to compare the dynamics to a model,
integrable, free dynamics.

Write 
$$\Phi :  \bR_t \times P \dashedrightarrow P $$
for the flow defined by our $n$-body system.  {\it We have used the
  broken arrow notation } for the map $\Phi$ to indicate that the
domain of the map need not be all of $ \bR_t\times P$, thus allowing
for the incomplete flows like the flows that occur for potentials such
as Newton's or Coulomb's which have singularities.  The curve
$t \mapsto \Phi_t (x)$, where defined, is a solution to our Hamilton's
equations having initial condition $x \in P$.
 
The {\em free flow} $\Phi^{(0)}$, on the other hand, is the flow whose
projected curves are the lines $t\mapsto a t + c$\,: 
\beq \Phi^{(0)}: \bR_t\times
P\to P\qmbox{,} \Phi^{(0)}_t(p, q)=(p,q+t\cM^{-1}p)   \Leq{free:solu}
and is generated by the free Hamiltonian $H_0 = K$. Let
\beq F_0 = F^+_0 := \big\{ (p,q)\in P \mid v= \cM^{-1}p\notin \Delta
\, \big\}. 
\Leq{P0free}

\begin{theorem} 
  \label{thm:mainMoller}\quad\\%
  Let $V$ be a short-range $(\alpha, k)$-potential with $\alpha>1$ and
  $k \ge 2$ (see definition~\ref{def:alpha-k}).\\
  Then the dynamics $\Phi$ on $F^+$ is 
  conjugate to the free dynamics $\Phi^{(0)}$: there exists a $C^{k-1}$
  symplectomorphism $\Omega : F_0 \to F^+$ %(onto its image, which need not be all of $F^+$) 
  such that
  $$\Omega \circ \Phi^{(0)}_t = \Phi_t \circ \Omega \qquad (\forall t
  \ge  0).$$
\end{theorem}
This is the qualitative contents of Theorem \ref{thm:moeller} below.
\medskip

An analogous theorem to \ref{thm:main2} holds for long-range
potentials.  Instead of comparing the given flow with the free flow,
we must compare it with an integrable, time dependent ``Dollard
Hamiltonian'' $H_D (p,t ) = K(p) + V((\sqrt{1+t^2})p)$ (which does not
depend on $q$!). See Theorem \ref{thm:Dollard:Moeller} for precise
statements.

\begin{earlier}\quad\\
  In 1927 {\sc Chazy} (\cite[Chapter 5]{Cha}) used the term
  ``hyperbolic'' in the classification the long-time behaviour of
  solutions in the long range case of celestial mechanics.  He
  established an analytic asymptotic expansion near infinity for his
  hyperbolic solutions with initial terms
  \beq
  q(t) = at + b \log(t) + c +\cO\big(\log(t)/ t\big) ; \quad b = +\nabla V(a)\quad
  \text{ as } t \to \infty 
  \Leq{hyperbolicChaz}
  Later, {\sc Saari} \cite[section 8]{Sa}, and   
  {\sc Marchal} and {\sc Saari} \cite[section 10]{MS}    extended
  and clarified Chazy's results, focussing on how 
  cluster energies and angular momenta approach their limits. Here
  ``cluster'' refers to the situation where $v^+ \in \Delta$. The 
  ``clusters'' are the subsets of mass indices $i$,
  for which $v^+ _i = v^+ _ j$.
  
  The $\log(t)$ term in Chazy's expansion
  equation~\eqref{hyperbolicChaz} is an essential consequence of the 
  $1/r$-nature of the Newtonian (or Coulomb) potential. On the other
  hand, hyperbolic solutions  for short range potentials   satisfy  
  \beq
  q(t) = a t + c + o(1), \quad\text{ as }t \to + \infty.
  \Leq{hyperShortRange}
  \textsc{Simon} \cite{Sim} proved the validity of this expansion for the
  two-body problem using the M\o ller transform (as used in
  section~\ref{sect:moeller}), or, as he called it, the
  {\em wave transformation}. 

  Recently {\sc Maderna} and {\sc Venturelli} \cite{MV} investigated
  forward hyperbolic motions for $n$-body problem using variational
  and weak KAM methods. And \textsc{Duignan} et al. set up
  \cite{Duignan} an approach to hyperbolic motions and scattering for
  the $n$-body problem which relies on a McGehee-style compactification
  of phase space which adds fixed points at infinity whose stable
  manifold correspond to forward hyperbolic solutions.\hfill $\Diamond$
\end{earlier}

\begin{subsection}{Summary: Main notations}
\label{subsec:notations} 

\begin{tabular}[t]{ll}
  $v^\pm(x)$ &asymptotic velocity of state $x$
             (definition~\ref{def:asymptotic:velocity})\\
  $P$ & phase space (equation~\eqref{phase_space})    \\     
  $\widehat P$ & phase space when collision singularities present (equation~\eqref{phase_space} )   \\     
  $F^\pm$ &forward and backward free regions
            (definition~\ref{def:free})\\
  $\widehat F^\pm$ &  as above,  but when collision singularities present
            (definition~\eqref{P:hat:free:pm})\\  
  $F^+_{\rm loc}$ &forward finally free region
                 (definition~\ref{def:finally:free})\\  
   $F_0$ &free region of the free flow (equation~\eqref{P0free})\\
  $\Phi^{(0)}_t$ & free flow (equation~\eqref{free:solu})\\
  $\Phi_t$ & n-body flow (subsection~\ref{subsec:flows})\\
  $\hat \Phi_t $ &n-body flow when collisions present \\
\end{tabular}

\end{subsection}

% ----------------------------------------------------
\section{Do we know when we are free?}
\label{sect:free} 
%----------------------------------------------------

For simplicity, we first consider long-range  potentials $V$ which are
non-singular at the origin i.e. $C^2$ on $\bR^{dn}$. Many
properties which hold for these potentials also hold for 
singular long range potentials (e.g.\ the
gravitational $n$-body potential). This will be proved in
section~\ref{sect:moeller2}.

We equip the real vector space of long range $\alpha$-potentials $V\in
C^2(\bR^{dn},\bR)$ (as introduced in definition~\ref{def:range},
with $\alpha>0$) with the seminorm  
\beq
\|V\|^{(\alpha)} \; := \;  \| \cM^{-1}\|\
\max_{i<j\in N}  \sup_{q\in\bR^d\backslash \{0\}} \|q\|^{\alpha+1}
\|\nabla V_{i,j}(q)\|.
\Leq{V:alpha} 
Typically, pair potentials are $C^2$, smoother, often even
analytic. In order to describe a section of the flow in restriction to
$F^+$, we will need a $C^2$ seminorm estimate of the potential
and later, scattering estimates will be improved using 
$C^k$-seminorms, with $k\in \bN$, $k \geq 2$. We now introduce such
seminorms.

\begin{definition}
  \label{def:alpha-k}\quad\\%
  An $(\alpha,k)$--\emph{potential} $V$ is a potential whose pair
  potentials $V_{i,j}\in C^k(\bR^d,\bR)$ fulfill
  \beq
  \partial^\gamma V_{i,j}(q)=\cO\big(\|q\|^{-\alpha - |\gamma|}\big)
  \qquad (\gamma\in\bN_0^d, |\gamma|\le k).
  \Leq{V:k:decay}
  On the space of $(\alpha,k)$-potentials we define
  \beq
  \|V\|^{(\alpha,k)} := \| \cM^{-1}\|\,
  \sum_{i<j\in N}  \sum_{\gamma\in\bN_0^d, |\gamma| = k}\
  \sup_{q\in\bR^d\backslash \{0\}} \|q\|^{\alpha+k} |\partial^\gamma
  V_{i,j}(q)|
  \Leq{V:k:norm}
  (so that $\|V\|^{(\alpha)} = \|V\|^{(\alpha,1)}$).
\end{definition}
The inessential factor $\|\cM^{-1}\|=m_{\min}^{-1}$ for
$m_{\min}:=\min(m_1,\ldots,m_n)$ simplifies formulae. 

Recall the definition of the free region
\beq
F^+= \big\{x\in P \mid v^+ (x) \in \bR^{dn}\, \, \backslash\, \Delta
\big\}. 
\Leq{P:free:pm:B}
It depends on the details of the (generally non-integrable)
flow, and hence implicitly on the precise form of the potential $V$.  
So general properties of the free region are hard to
grasp. Surprisingly, there is an explicit surface of  
section of the flow restricted to $F^+$. It bounds a
positive-invariant subset $F^+_{\rm loc} \subset F^+$ having the
property that every orbit in $F^+$ must eventually enter $F^+_{\rm loc}$. 

We still assume that $V$ is a non-singular $\alpha$-potential. Let
$\delta$ and $C$ with 
\[0<\delta \le \delta_0:= \min(\alpha/(4+\alpha),1/5), \quad C := 16
\,dn\,\|V\|^{(\alpha{,2})}/\delta.\]
Define $F^+_{\rm loc}$ by \eqref{def:ff}, with our chosen values of
$\alpha$, $\delta$ and $C$. The three inequalities assert
\begin{enumerate}[--]
\item the dominance of the interparticle kinetic energy over potential
  energy
\item the  near-parallelism of interparticle distances and velocities
\item that interparticle  distances are nearly proportional to
  interparticle  velocities. 
\end{enumerate}
This tells us what the motion of free particles eventually looks like.
For example, landing in $F^+_{\rm loc}$ yields the simple propagation estimates
\eqref{propagation:est}. 

\begin{theorem}[Final free region]
  \label{thm:final:free}%
  \quad
  \begin{enumerate}[1.]
  \item $F^+_{\rm loc}$ is forward invariant : $\Phi_t (F^+_{\rm loc}) \subseteq
    F^+_{\rm loc}$ for $t\ge0$.
  \item $F^+_{\rm loc}$ is a subset of $F^+$.
  \item For any $x_0\in F^+$ there is a time $t$ such that
    $\Phi(t,x_0)\in F^+_{\rm loc}$.
  \item 
    For $x_0\in F^+_{\rm loc}$ the distance between the particles
    $i<j\in N$ increases linearly:
    \beq
    \eh v_{i,j}(0) \; t\le q_{i,j}(t,x_0) - q_{i,j}(0,x_0) 
    \le  \textstyle\frac{3}{2} v_{i,j}(0)\; t
    \qquad\big(t\in[0,\infty)\big).
    \Leq{propagation:est} 
  \end{enumerate}
\end{theorem}

As already mentioned, the boundary of $F^+_{\rm loc}$ is thus a $C^0$
surface of section of the flow restricted to the free region.

\begin{example}[Two bodies]
  We already remarked that $F^+_{\rm loc}$ is invariant under 
  Euclidean transformations. So in particular 
  $F^+_{\rm loc} = T^*\bR^d\times \widetilde{F}^+_{\rm loc} $, the cartesian product referring to 
  the separation of center of mass and internal motion.\\
  In the case $n=2$,  $\widetilde{F}^+_{\rm loc}$
  is fibered over $\bR^{d}_q$, with fiber
  diffeomorphic to the cylinder $[1,\infty)\times B_{d-1}$ in spherical coordinates, 
  see figure \ref{fig:F:plus}.
  \hfill $\Diamond$
  
  \begin{figure}[h]
    \begin{center}
      \includegraphics[width=8cm]{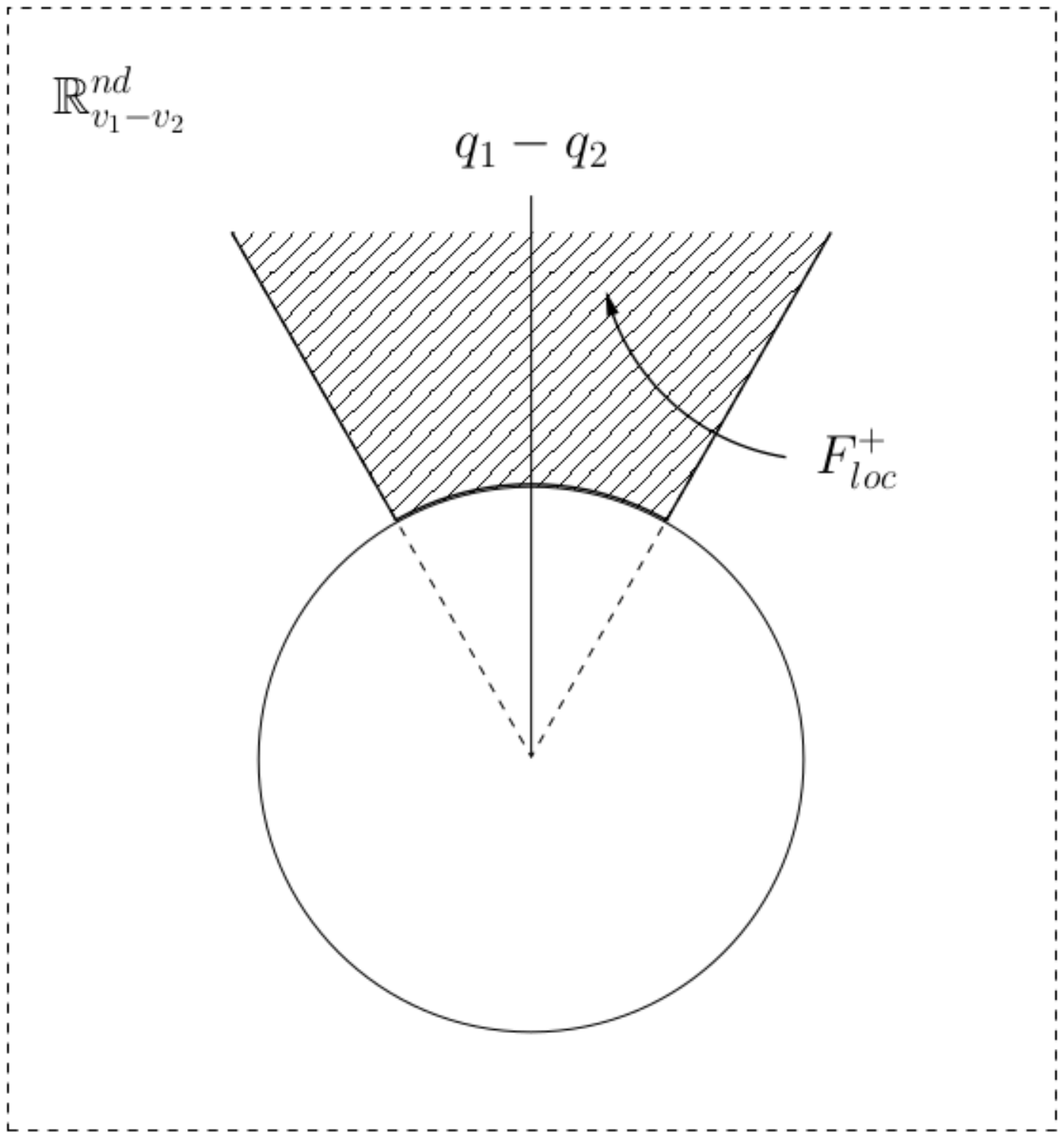}
    \end{center}
    \caption{Fiber over $(q,v_1)$ of $\partial F^+_{\rm loc}$ in the $n=2$
    case}
    \label{fig:F:plus}
  \end{figure}
\end{example}
\begin{remark}[Topology of the final free region]\quad\\ 
Although $\Delta$ is contractible,
already the set $\bR^{dn}\, \backslash\, \Delta$ to which partial $F^+_{\rm loc}$ projects, 
is topologically rich:
\begin{enumerate}[1.]
\item 
For $d=1$, there is a homeomorphism
$\bR^n \backslash\, \Delta \cong \bR^n\times {\rm Sym}(n)$.
\item 
For $d=2$, the cohomology ring of $\bR^{2n}\, \backslash\, \Delta$ is 
the one of the product over $k$ of bouquets of  $1\le k\le n-1$ circles,
see {\sc Arnold} \cite{Ar}.
\hfill $\Diamond$
\end{enumerate}
\end{remark}
\medskip 
In order to prove Theorem \ref{thm:final:free}, we will use the following
lemma, whose proof is routine, and where we denote by $\LA\, \cdot\, \RA$ a 
smoothened version a the absolute value: 
\[\LA q \RA =\sqrt{q^2+1}.\]
\begin{lemma}
  \label{lem:integral:estimate}%
  \quad\\[-6mm]
  \begin{enumerate}[1.]
  \item 
    For $\alpha>0$ and $q\ge 0$,
    \beq 
    {\textstyle \frac{1}{\alpha}}\LA q \RA^{-\alpha} 
    \ \le \ \int_{q}^{\infty}\! \LA \tilde{q} \RA^{-\alpha-1}d\tilde{q} 
    \ \le \ (\textstyle \frac{1}{\alpha}+1)\LA q \RA^{-\alpha}.
    \Leq{ineq:smooth:abs} 
  \item If $\,V$ is an $\alpha$-potential (see~\eqref{V:decay}),
    \[\|V\|^{(\alpha)} \le \|V\|^{(\alpha,1)}\]
    and, if $\,V$ is an $(\alpha,k+1)$--potential
    (see~\eqref{V:k:decay}) with $k\geq 1$, 
    \[ \|V\|^{(\alpha,k)} \le  d\, \|V\|^{(\alpha,k+1)}.\]
  \end{enumerate}
\end{lemma}

\noindent
\textbf{Proof of Theorem \ref{thm:final:free}:} \\
We will repeatedly use   the symbol $X_{i,j}$ for the  relative
accelerations
\[X_{i,j}: \bR^{dn}\to \bR^d \qmbox{,}
X_{i,j}(q):=\sum_{k\in N\setminus\{i\}}\!\!\!\frac{\nabla V_{i,k}(q_i-q_k)}{m_i}-
\sum_{k\in N\setminus\{j\}} \!\!\! \frac{\nabla V_{j,k}(q_j-q_k)}{m_j}\]
between  the $i$-th and $j$-th particle, and    the estimate
\beq
 \|X_{i,j}(q)\| \le 2(n-1)\frac{\|V\|^{(\alpha)}}{ q_{\min}^{\alpha+1}}\qquad(i<j\in N).
\Leq{X:est}  
Throughout the proof we also use that, by Lemma
\ref{lem:integral:estimate}.2,  
\beq
C\ge 8 n\|V\|^{(\alpha)}/\delta.
\Leq{C:ge}
\begin{enumerate}[1.]
\item $F^+_{\rm loc}$ is open, since it is defined by strict inequalities
  among continuous functions on phase space. To prove that $F^+_{\rm loc} $
  is forward invariant, it is sufficient to show that the Hamiltonian
  vector field {\em points inwards along its boundary}
  $\partial F^+_{\rm loc}$. Thus we will show that the difference of the
  sides of each inequality has positive time derivative (Poisson
  bracket $\{\bullet,H\}$) at instants at which that inequality
  becomes an equality.

  Note that on $\overline{F^+_{\rm loc}}$ the phase space functions
  $q_{i,j}$ and $v_{i,j}$ have positive values and are thus smooth.
  For $q_{\min}$, $q_{\max}$, $v_{\min}$ and $v_{\max}$, which are
  only Lipschitz continuous, we consider the distributional
  derivative.
\begin{enumerate}[(a)]
\item On $F^+_{\rm loc}$ the time derivative relating to the first
  inequality in \eqref{def:ff} is positive, since defining $\{f, H\}$,
  as the distribution $\frac{d}{dt}f$ along orbits, $v_{\min}$
  fulfills the inequality
  \[\{v_{\min}^2, H\}\ge -2v_{\min}\, \max_{i<j\in
    N}\|\dot{v}_{i}-\dot{v}_{j}\|,\]
  see {\em e.g.}\ {\sc Lieb} and {\sc Loss} \cite[Cor. 6.18]{LL} for
  the weak gradient of the minimum of functions. Thus
\begin{align}
\Big\{v_{\min}^2 -&\, C \frac{q_{\max}}{q_{\min}^{\alpha+1}} \,,\, H\Big\}\ge\NN\\
\ge&\; C(\alpha+1) \frac{q_{\max}}{q_{\min}^{\alpha+2}}(1-\delta)v_{\min}
-2v_{\min}\, \max_{i<j\in N}\|X_{i,j}(q)\| -
C \frac{v_{\max}}{q_{\min}^{\alpha+1}}\NN\\
\ge&\;  C\Big((\alpha+1)(1-\delta)\frac{v_{\min}}{q_{\min}} -\delta/4  
\frac{v_{\min}} {q_{\max}}- \frac{v_{\max}} {q_{\max}}\Big)
\frac{q_{\max}}{q_{\min}^{\alpha+1}}\NN\\
\ge&\;  C\Big((\alpha+1)(1-\delta)-\delta/4 -(1+2\delta)\Big)
\frac{v_{\min}\,q_{\max}}{q_{\min}^{\alpha+2}}\NN\\
=&\;  C\Big(\alpha(1-\delta) -\frac{13}{4}\delta\Big)
\frac{v_{\min}\,q_{\max}}{q_{\min}^{\alpha+2}}> 0.\NN
\end{align}
The factor $1-\delta$ in the first inequality follows from the second
line of \eqref{def:ff}. 
The second inequality follows from \eqref{X:est} and \eqref{C:ge}.
The factor $1+2\delta$ in the second to last line follows from the
third line of \eqref{def:ff}.\\ 
In the final inequality we used that $\delta \in
\big(0,\min(\alpha/(4+\alpha),1/5)\big]$: 
\begin{enumerate}[$\bullet$]
\item 
For $\alpha\in (0,1]$ we obtain $\alpha(1-\delta) -\frac{13}{4}\delta \ge \frac{3\alpha}{4(4+\alpha)}>0$.
\item 
For $\alpha\in (1,\infty]$ we get 
$\alpha(1-\delta) -\frac{13}{4}\delta\ge\frac{4}{5}\alpha- \frac{13}{20}\ge\frac{3}{20} $. 
\end{enumerate}
\item 
The time derivative of the left hand side of the second inequality
\[ \LA v_i-v_j, q_i-q_j\RA - (1-\delta)v_{i,j}q_{i,j} > 0 \] 
in \eqref{def:ff} is positive, too. This is trivial if
$\|V\|^{(\alpha)}=0$, that is, $V=0$. Otherwise 
\begin{align}
\big\{\langle v_i-v_j,&\, q_i-q_j\rangle  - (1-\delta)v_{i,j}q_{i,j} \;,\; H\big\}\ge \NN\\
\ge&\; v_{i,j}^2 - \LA X_{i,j}(q),q_i-q_j\RA  - (1-\delta)\big(v_{i,j}^2+\|X_{i,j}(q)\|q_{i,j}\big)
\NN\\
\ge&\; \delta C \frac{q_{\max}}{q_{\min}^{\alpha+1}}-
2(n-1)\frac{\|V\|^{(\alpha)}q_{\max}}{ q_{\min}^{\alpha+1}}
(2-\delta)
\NN\\
\ge &\; \big(8n-2(n-1)(2-\delta)\big)\frac{\|V\|^{(\alpha)}q_{\max}}{ q_{\min}^{\alpha+1}}
> 4n\frac{\|V\|^{(\alpha)}q_{\max}}{ q_{\min}^{\alpha+1}}
 >0 .\NN
\end{align}
For the third line we used the first inequality of \eqref{def:ff} and
\eqref{X:est}, and  
\eqref{C:ge} for the last line.
\item 
Concerning the third  inequality in \eqref{def:ff},
as $\{ q_{k,l} , H \}\in [(1-\delta)v_{k,l}, v_{k,l}]$ (using the
second inequality in \eqref{def:ff}),  
at value zero of 
$(1+2\delta) \frac{q_{k,l}}{v_{k,l}}- \frac{q_{i,j}}{v_{i,j}}$ its
time derivative estimates as follows: 
\begin{align}
&\Big\{ (1+2\delta) \, \frac{q_{k,l}}{v_{k,l}}- \frac{q_{i,j}}{v_{i,j}}, H  \Big\}\ge\NN\\
\ge&\;  (1+2\delta)(1-\delta) -1 -( (1+2\delta) q_{k,l}\|X_{k,l}(q)\| +q_{i,j}\|X_{i,j}(q)\|)/v_{\min}^2\NN\\
\ge&\;  \delta-2\delta^2- 
2 (n-1) (1+\delta) \frac{\|V\|^{(\alpha)}q_{\max}}{v_{\min}^2 \, q_{\min}^{\alpha+1}}
\NN\\
>&\textstyle  \; \delta(1- 2\delta) - 4n(1+\delta) \frac{\|V\|^{(\alpha)}}{C  }
\ge \delta(1- 2\delta) -\eh(1+\delta))  
\ge\delta(\frac{3}{5}\! -\! \frac{6}{10}) = 0. \NN
\end{align}
\end{enumerate}
\item[4.]
  We now prove item 4, before items 2 and 3. Throughout the proof of
  item 4 we will use that we already proved positive invariance of
  $F^+_{\rm loc}$ (item 1). We adopt the notation
  $\tilde{f}(t):=f\circ\Phi(t,x_0)$ for a phase space function $f$,
  with $x_0\in F^+_{\rm loc}$.

  For $F:=\eh \tilde{q}_{i,j}^{\,2}$ and $t\ge0$ we get from
  \eqref{def:ff} that \beq F'(t)\ge
  (1-\delta)\tilde{q}_{i,j}(t)\tilde{v}_{i,j}(t) \Leq{Fp} and \beq
  F''(t)\ge (1-\delta)
  \l[\tilde{v}_{i,j}^{\,2}(t)-\tilde{q}_{i,j}(t)\|
  \tilde{X}_{i,j}(t)\|\ri]
  \ge(1-\delta)(1-\delta/4)\tilde{v}_{i,j}^{\,2}(t) \ge \textstyle
  \frac{19}{25} \tilde{v}_{i,j}^{\,2}(t).  \Leq{Fpp} The second
  inequality in \eqref{Fpp} is valid, since by \eqref{X:est},
  \eqref{C:ge} and \eqref{def:ff} \beq \tilde{q}_{i,j}(t) \|
  \tilde{X}_{i,j}(t) \| \ \le\ \tilde{q}_{i,j}(t)
  2(n-1)\frac{\|V\|^{(\alpha)}}{ \tilde{q}_{\min}^{\alpha+1}(t)} \
  \le\ \ev C\frac{\tilde{q}_{i,j}(t)}{\tilde{q}_{\min}^{\alpha+1}(t)}
  \delta \ \le\ \frac{\delta}{4}\tilde{v}_{\min}^2(t).  \Leq{tq:X} The
  third inequality in \eqref{Fpp} follows from $\delta \le1/5$.

There exists a maximal $T\in (0,+\infty]$ so that
\beq
\big\|\big(\tilde{v}_i(t)-\tilde{v}_j(t)\big) - \big(\tilde{v}_i(0)-\tilde{v}_j(0)\big)\big\|^2
\le \textstyle\frac{1}{6} \tilde{v}_{i,j}^2(0)\qquad \big( t\in [0,T) \big).
\Leq{assum:velocity} 
Thus $\big(1-\sqrt{1/6}\big)^2\,\tilde{v}_{i,j}^{\,2}(0)\le \tilde{v}_{i,j}^{\,2}(t)\le 
\big(1+\sqrt{1/6}\big)^2\,\tilde{v}_{i,j}^{\,2}(0)$
within this time interval, and by \eqref{Fp}, \eqref{Fpp} this implies
\beqn
F(t)&=&\textstyle F(0)+\int_0^t \big(F'(0)+\int_0^sF''(\tau)d\tau\big)\, ds
\label{F:FFF}\\
&\ge& F(0) +(1-\delta)\tilde{q}_{i,j}(0)\tilde{v}_{i,j}(0) t 
+\textstyle\frac{19}{25} \int_0^t\int_0^s \tilde{v}_{i,j}^{\,2}(\tau)
\,d\tau\,ds\NN\\ 
&\ge& \textstyle
\eh \tilde{q}_{i,j}^{\,2}(0)+\frac{4}{5}\tilde{q}_{i,j}(0)\tilde{v}_{i,j}(0)t
+\frac{19}{50} \big(1-\sqrt{1/6} \big)^2\,\tilde{v}_{i,j}(0)^2 t^2. \NN
\eeqn
Conversely by the first line in \eqref{F:FFF} and \eqref{tq:X} with
$\delta\le 1/5$, 
\beqn
F(t)&\le& F(0) +\tilde{q}_{i,j}(0)\tilde{v}_{i,j}(0) t
+\textstyle\frac{21}{20}  
\int_0^t\int_0^s \tilde{v}_{i,j}^{\,2}(\tau) \,d\tau\,ds  \NN \\
&\le& F(0) +\tilde{q}_{i,j}(0)\tilde{v}_{i,j}(0) t +\textstyle 
\frac{21}{40} \big(1+\sqrt{1/6}\big)^2\,\tilde{v}_{i,j}^{\,2}(0)t^2. \NN
\eeqn
These two estimates prove both inequalities in 
\eqref{propagation:est} for time $t\in[0,T)$.

\item[2.] Next we start by showing that $T =+\infty$ in
  \eqref{assum:velocity}.  With the rescaled time parameter
  \[s(t) := \textstyle\frac{\tilde{v}_{i,j}(0)} {2\tilde{q}_{i,j}(0)}
  t\qmbox{,} \tilde{q}_{i,j}(t)\ge \tilde{q}_{i,j}(0)\LA s(t)\RA.\]
  Note that by definition \eqref{def:ff} of $F^+_{\rm loc}$ the scaling
  factors $\textstyle\frac{\tilde{v}_{i,j}(0)} {2\tilde{q}_{i,j}(0)}$
  are, up to a factor $1+\delta$, independent of the index pair
  $(i,j)$. So by applying \eqref{assum:velocity}, \eqref{X:est}, and
  \eqref{def:ff} with \eqref{C:ge} in succession, 
  \beqno \lefteqn{
    \big\|\big(\tilde{v}_i(t)-\tilde{v}_j(t)\big) - \big(\tilde{v}_i(0)-\tilde{v}_j(0)\big)\big\|^2=}\\
  &=&-2\int_0^t\LA \tilde{v}_{i}(\tau)
  -\tilde{v}_{j}(\tau) - \big(\tilde{v}_i(0)-\tilde{v}_j(0)\big),\tilde{X}_{i,j}(\tau)\RA\,d\tau\\
  &\le& \textstyle
  \frac{2}{\sqrt{6}}\,\tilde{v}_{i,j}(0)\int_0^\infty \|\tilde{X}_{i,j}(\tau)\|\,d\tau\\
  &\le& \textstyle \frac{4}{\sqrt{6}}\,(n-1)\,\tilde{v}_{i,j}(0)
  \|V\|^{(\alpha)}
  \int_0^\infty  \tilde{q}_{\min}(\tau)^{-\alpha-1}\, d\tau\\
  &\le& \textstyle \frac{8}{\sqrt{6}}\, (n-1)\,
  \tilde{v}_{i,j}(0)\|V\|^{(\alpha)} \frac{(1+\delta)
    \tilde{q}_{\max}(0)}
  {\tilde{q}_{\min}^{\alpha+1}(0)\tilde{v}_{i,j}(0)}
  \int_0^\infty  \LA s\RA^{-\alpha-1}\, ds\\
  &\le& \textstyle \frac{\sqrt{6}}{5}\,\tilde{v}_{i,j}^{\,2}(0) \,
  \delta
  \int_0^\infty  \LA s\RA^{-\alpha-1}\, ds\\
  &\le& \textstyle \frac{\sqrt{6}}{5} \tilde{v}_{i,j}^{\,2}(0)
  \min\!\big(1/5,\frac{\alpha}{4+\alpha}\big)
  \frac{\sqrt{\pi} \Gamma(\alpha/2)} {2\Gamma((1+\alpha)/2)}\\
  &\le& \textstyle \frac{\sqrt{6}\pi}{50}\tilde{v}_{i,j}^{\,2}(0) <
  \tilde{v}_{i,j}(0)^2/6\, , \eeqno since
  $\min\!\big(1/5,\frac{\alpha}{4+\alpha}\big) \frac{\sqrt{\pi}
    \Gamma(\alpha/2)} {2\Gamma((1+\alpha)/2)}$
  attains its maximal value $\pi/10$ for $\alpha=1$.  This shows that
  in \eqref{assum:velocity} $T=+\infty$.  Thus by
  \eqref{assum:velocity} the velocity differences stay bounded away
  from zero ($\tilde{v}_{i,j}(t)\ge \eh \tilde{v}_{i,j}(0)>0$ for all
  $t\ge0$) so that the initial condition $x_0\in F^+_{\rm loc}$ is in
  $F^+$.
\item[3.] Let $x_0\in F^+$. By definition, $v^+(x_0)$ exists and
  $v^+ (x_0) \notin \Delta$.  It follows that for any
  $\delta\in(0,\delta_0]$ there exists a time $t_0$ such that \beq
  \|\tilde{v}_k(t)-v^+_k(x_0)\| \le \ea\,\delta \,
  \ov{v}^+_{\min}(x_0) \qquad (k\in N,\,t\ge t_0).  \Leq{v:ofv:diff}
  In particular
  \[\tilde{v}_{\min}(t)\ge (1-\ev\delta)\,  \ov{v}^+_{\min}(x_0) > 0 \qquad (t\ge t_0).\]
  As
  $ \big\|\big( \tilde{q}_{i}(t)-\tilde{q}_{j}(t)\big) -\int_{t_0}^t
  (\tilde{v}_i(s)-\tilde{v}_j(s))\,ds\big\| = \tilde{q}_{i,j}(t_0)$,
  \beq \big\| \big( \tilde{q}_{i}(t)-\tilde{q}_{j}(t)\big)\, -\,
  (t-t_0)\big(\ov{v}_i(x_0)-\ov{v}_j(x_0)\big)\big\| \le \delta/4 \,
  (t-t_0) \ov{v}^+_{\min}(x_0) + \tilde{q}_{i,j}(t_0).
  \Leq{near:affine} So
  $v_{\min}(x)^2> C \frac{q_{\max}(x)}{q_{\min}(x)^{\alpha+1}}$ for
  $x:=\Phi(t,x_0)$, $t\ge t_0$ large. This is the first condition in
  the definition \eqref{def:ff} of $F^+_{\rm loc}$.

  Concerning the second condition, similarly by \eqref{near:affine},
  for $t$ large
  \[\LA \tilde{v}_{i}(t)-\tilde{v}_{j}(t),
  \tilde{q}_{i}(t)-\tilde{q}_{j}(t)\RA\ge
  (1-\delta)\tilde{v}_{i,j}(t)\tilde{q}_{i,j}(t)\]
  and, for the third condition,
  \[(1+2\delta) \frac{\tilde{q}_{k,l}(t)}{\tilde{v}_{k,l}(t)} >
  \frac{\tilde{q}_{i,j}(t)}{\tilde{v}_{i,j}(t)} \qquad(\,i<j,k<l \in
  N).\]
  This shows that $\Phi_t (x_0)\in F^+_{\rm loc}$ for all $t$ sufficiently
  large.  \hfill $\Box$%\\[2mm]
\end{enumerate}

%----------------------------------------------------
\section{Regularity of the asymptotic velocity}
%----------------------------------------------------

We move on to the
regularity of the asymptotic velocity map $v^+ : x \mapsto v^+ (x)$.

\begin{theorem}[{{\cite[Theorem 5.4.1]{DG}}}]
Let the potential $V \in C^2(\bR^{dn},\bR)$ be long range. Then the 
asymptotic velocity $v^+(x)$ exists for all $x \in P$.
\end{theorem}

The map $v^+ : P \to \bR^{dn}$ is Borel-measurable, but may be
discontinuous. 

\begin{example}[Discontinuity of the asymptotic velocity]
  \label{rem:discontinuous:as}%
  \quad\\
  Take $d=1$ and $n=2$ and a non-negative pair potential
  $V_{1,2}\in C^2_{\rm c}(\bR,\bR)$ which has compact support, and a
  unique maximum $V(0)>0$. The velocity along any trajectory is
  constant after some time. So, the map $v^+$ is defined on the whole
  phase space $\bR^2$ and has a discontinuity at the hyperbolic
  equilibrium $x=0$ and nowhere else.
  \hfill $\Diamond$
\end{example}

We will see that in restriction to the free region $F^+$ the map
$v^+$ is continuous and even differentiable. We will use the notation
\[p^+(x_0)  = \cM v^+(x_0)= \lim_{t\rightarrow\infty}p(t).\]

\begin{theorem}
  \label{thm:both:moeller}
  \quad\\ 
  Let $V$ be an $(\alpha,k)$--potential.  
  \begin{enumerate}[1.]
  \item \label{Statement1} The map $x_0 \mapsto v^+ (x_0)$ is a
    $C^{k-1}$ map $F^+ \to \bR^{dn}$.
  \item \label{Statement2}
    Quantitatively, if $x_0 = (p_0,q_0) \in F^+_{\rm loc}$, for 
    multi-indices $\delta:=(\beta,\gamma)\in\bN_0^{dn}\times \bN_0^{dn}$ 
  with $|\delta|\equiv |\beta|+|\gamma| \le k - 1$ and
  partial derivatives:  $\pa^\delta_{x_0}:=\pa^\beta_{p_0}\pa^\gamma_{q_0}$ we get
    \beq
    \pa^\delta_{x_0} (p^{+} (x_0) - p_0) =
    \cO\l(\|V\|^{(\alpha,k)} v_{\min}(x_0)^{-1-|\beta|}  \LA
    q_{\min}(x_0)\RA^{-\alpha-|\gamma|} \ri). 
    \Leq{O:p}
  \end{enumerate}
\end{theorem}

\begin{remark}[Variants]\quad\\
  The constant in the order estimate \eqref{O:p} is independent of
  $V$.  Using the first condition in the definition \eqref{def:ff} of
  $F^+_{\rm loc} $, we obtain the weaker estimate
  \[\pa^\delta_{x_0} (p^{+} (x_0) - p_0) =
  \cO\l( v_{\min}(x_0)^{+1-|\beta|} \LA q_{\min}(x_0)\RA^{-|\gamma|}
  \ri).\]
  Similarly, instead of \eqref{Moe:E}, we would have the weaker
  estimate
  \[\pa^\delta_{X_0} (Q_0-q_0) =
  \cO\l( v_{\min}(X_0)^{-|\beta|} \LA q_{\min}(X_0)\RA^{1-|\gamma|}
  \ri).\]
  These estimates depend on the norm $\|V\|^{(\alpha,k)}$ of the
  potential only indirectly, via the phase space region $F^+_{\rm loc} $
  where they apply.\hfill $\Diamond$
\end{remark}
\textbf{Proof of Theorem \ref{thm:both:moeller}:}\\
We use the shorthands $q_{\min}:= q_{\min}(x_0)$, $v_{\min}:= v_{\min}(x_0)$
and continue to use the notation 
$\tilde{f}(t):=f\circ\Phi(t,x_0)$ for a phase space function $f$.\\
$\bullet$
To prepare for the proof of Claim \ref{Statement1}, 
we first estimate the initial value problem for long-range potentials.
As $V$ is an $(\alpha,k)$--potential, the flow 
\[\Phi\in C^{k-1}(\bR\times P,P).\]
For derivatives $\pa^\delta_{x_0}$ w.r.t.\ initial conditions $x_0$ with $1\le |\delta|\le k-1$, 
like in \cite[section 6]{Kn} we use the integral representation of the trajectory
\[ q(t,x_0)= q_0+
\cM^{-1}\l(\textstyle t p_0 - \int_0^t\!\! \int_0^s \nabla V\big(q(\tau,x_0)\big)\,d\tau \,ds \ri)
\qquad \big(t\in[0,\infty)\big).\]
By a standard dominated convergence argument (see, {\em e.g.}, 
{\sc Elstrodt} \cite[Thm.\ IV.5.7]{El}) its deviation from free motion is controlled by
\beqn
\lefteqn{\pa^\delta_{x_0}
\big( q(t,x_0) -(q_0+t  \cM^{-1}p_0)\big)= -\int_0^t \!\! \int_0^s \!
 \cM^{-1} \pa^\delta_{x_0}\nabla V\big( q(\tau,x_0)\big)\,d\tau \,ds 
=} &&\label{implicit:equation}\\
&&\hspace*{-8mm} - \! \sum_{N=1}^{|\delta|} \label{eq:schlange}
 \cM^{-1}\hspace*{-7mm}
 \sum_{\stackrel{\delta^{(1)}+\ldots+\delta^{(N)} = \delta} {|\delta^{(i)}|>0}}\!\!\!\!\!
\int_0^t\!\! \int_0^s \!
D^N\nabla V\big(q(\tau,x_0)\big)
\l(\pa^{\delta^{(1)}}_{x_0}q(\tau,x_0),\ldots, 
\pa^{\delta^{(N)}}_{x_0}q(\tau,x_0)\ri) d\tau \,ds .\NN
\eeqn
Due to the $N=1$ term this is only an implicit equation for $\pa^\delta_{x_0}q(t,x_0)$.
To transform it into an explicit equation, we 
thus consider for  $\lambda>0$ the real Banach space 
$\big(\widehat{\cC}, \|\cdot\|_{\lambda}\big)$,
\beq
\widehat{\cC} := \l\{ w\in C\l([0,\infty),\bR^{dn}\ri) \l| \,
\|w\|_{\lambda}:=\sup_{t\geq0} \| w(t)\|/\langle \lambda  t\rangle < 
\infty \ri.\ri\},
\Leq{hat:C}
noting that $\widehat{\cC}$ is independent of the choice of $\lambda$.
The linear operator ${\cal Q}\equiv {\cal Q}_{x_0}$, 
\beq
{\cal Q}(w)(t):=
\cM^{-1}  \int_0^t \!\! \int_0^s  D\nabla V\big(q(\tau,x_0)\big)w(\tau)\,d\tau \,ds
\qquad(t\geq 0),
\Leq{Q:op}
maps  $\widehat{\cC}$ into itself, 
and we want to prove that for all $x_0\in F^+_{\rm loc}$
the operator norm of ${\cal Q}_{x_0}$ is strictly smaller than one for a 
suitable $\lambda$. Using \eqref{propagation:est} and \eqref{V:k:norm}, 
the operator norm is estimated by
\beqno
\|{\cal Q}\|_{\lambda}\!\!
&:=& \sup_{w:\,\|w\|_{\lambda}=1} \! \! \|{\cal Q}(w)\|_{\lambda}
\le {\|V\|^{(\alpha,2)}}
\sup_{t\ge 0} \frac{ \int_0^t \!\! \int_0^s 
( q_{\min}+\eh v_{\min}\tau )^{-2-\alpha}
\langle \lambda  \tau\rangle\,d\tau \,ds}{\langle \lambda  t\rangle}\\
&\le& \frac{\|V\|^{(\alpha,2)}}{\lambda^2  q_{\min}^{2+\alpha}}
\sup_{t\ge 0} \frac{ \int_0^t \!\! \int_0^s 
\langle \tau\rangle^{-1-\alpha}\,d\tau \,ds}{\langle t\rangle}
\le 4 (1+1/\alpha) \frac{\|V\|^{(\alpha,2)}}{ v_{\min}^2 q_{\min}^{\alpha}}
\eeqno
using Lemma \ref{lem:integral:estimate} in the last inequality and
setting 
\[\lambda:= \eh v_{\min}/q_{\min} .\] 
By Definition \eqref{def:ff} of $F^+_{\rm loc}$ the operator is a contraction:
\[\|{\cal Q}\|_{\lambda}\le \frac{4(1+1/\alpha)}{16\pi dn\,\max(1,1/\alpha)}\le 
\frac{1}{2\pi dn} < 1.\]
Thus \eqref{implicit:equation} can be transformed into
\beqn
\lefteqn{(\idty+ {\cal Q})(\pa^\delta_{x_0}q)(t) = \pa^\delta_{x_0}(q_0+t \cM^{-1}p_0)\ - 
\ \cM^{-1}\times}&&
\label{oneminq}\\
&&\hspace*{-6mm}\sum_{N=2}^{|\delta|}
\sum_{\stackrel{\delta^{(1)}+\ldots+\delta^{(N)}=\delta}{|\delta^{(i)}|>0}}\!\!\!\!\!
\int_0^t \!\! \int_0^s 
D^{N} \nabla V\big( q(\tau,x_0)\big)
\l(\pa^{\delta^{(1)}}_{x_0} q(\tau,x_0),\ldots, 
\pa^{\delta^{(N)}}_{x_0} q(\tau,x_0)\ri) d\tau \,ds\NN
\eeqn
with the  invertible operator $\idty + {\cal Q}$ on $\widehat{\cC}$.
As on the r.h.s.\ of (\ref{oneminq}) only partial derivatives 
of order $|\delta^{(i)}|< |\delta|$ appear, we can perform an 
induction in $|\delta|$.\\

Assume that for all $\delta'=(\beta',\gamma')\in\bN_0^{dn}\times \bN_0^{dn}$
with $1 \le |\delta'|\le |\delta|-1$
\beq
 \|\pa^{\delta'}_{x_0}q(\cdot,x_0) \|_\lambda =
\cO\l(  v_{\min}(x_0)^{-|\beta'|} \, q_{\min}(x_0)^{1-|\gamma'|}\ri).
\Leq{ass:delta:bar}
This assumption is satisfied for $|\delta'|=1$, since then the 
sum on the r.h.s.\ of \eqref{oneminq} equals zero. 
Then by \eqref{propagation:est} and \eqref{ass:delta:bar} the terms
on the r.h.s.\ of (\ref{oneminq}) fulfill  
\begin{align}
\MoveEqLeft 
\l\|\cM^{-1}  \int_0^t\!\! \int_0^s \! D^N \nabla V\big(q(\tau,x_0)\big)
\l(\pa^{\delta^{(1)}}_{x_0}q(\tau,x_0),\ldots, 
\pa^{\delta^{(N)}}_{x_0}q(\tau,x_0)\ri) d\tau \,ds\ri\| \le 
\label{so:wie}\\
&\le
 \|V\|^{(\alpha,N+1)} 
 \int_0^t\!\! \int_0^s \! q(\tau,x_0)^{-\alpha-N-1}
\prod_{i=1}^N \|\pa^{\delta^{(i)}}_{x_0}q(\tau,x_0)\| \,\, d\tau \,ds \NN\\
&\le
 \|V\|^{(\alpha,N+1)}\ \times \NN\\
&\ \int_0^t\!\! \int_0^s \! \big(q_{\min}(x_0)+\eh v_{\min}(x_0)t \big)^{-\alpha-N-1}
\prod_{i=1}^N \l(\|\pa^{\delta^{(i)}}_{x_0}q(\cdot,x_0)\|_\lambda \LA \lambda\tau \RA\ri)
\,\, d\tau \,ds \NN\\
&\le 
C_0 \|V\|^{(\alpha,N+1)} \, v_{\min}(x_0)^{-|\beta|} \, q_{\min}^{-\alpha-|\gamma|-1} 
\int_0^t\!\! \int_0^\infty \LA \lambda \tau\RA^{-\alpha-1} \, d\tau\,ds \ \NN\\
&\le  
C_1 \|V\|^{(\alpha,N+1)} \,
v_{\min}^{-2-|\beta|} \, q_{\min}^{1-\alpha-|\gamma|}\LA \lambda t\RA. \NN
\end{align}
For $x_0\in F^+_{\rm loc}$ that term is bounded above by 
(see \eqref{def:ff}) %
$C_\delta \, v_{\min}^{-|\beta|}   q_{\min}^{1-|\gamma|}\LA \lambda t\RA$,\footnote{$\delta$
of $C_\delta$ does not refer to the multi-index $\delta\in \bN_0^{2dn}$, but to the constant in
Theorem \ref{thm:final:free}. It is chosen as 
$\delta:=\min(\delta_0,\alpha-1)$ in the short range case ($\alpha>1$)
and $\delta:=\delta_0$ if $0<\alpha\le 1$.}
proving the induction step for \eqref{ass:delta:bar}.\\
$\bullet$
We prove the momentum estimate in \eqref{O:p} for no partial derivative
w.r.t.\ initial conditions ($\delta=0$), which holds for all long range potentials.
By the propagation estimate \eqref{propagation:est}
uniformly in $t\ge 0$
\beqn
\lefteqn{\|  \cM^{-1}(\tilde{p}(t)-\tilde{p}(0))\| } \NN\\
&\le& \int_0^t \|\cM^{-1}\nabla V(\tilde{q}(s))\|\,ds
\le \|V\|^{(\alpha,1)}\int_0^t \big( q_{\min} +  \eh v_{\min}\ s \big)^{-\alpha-1}\,ds \NN\\
&\le& 
\frac{\|V\|^{(\alpha,1)}}{\eh v_{\min}} \int_0^\infty \big( q_{\min}  +  s \big)^{-\alpha-1}\,ds 
=  \frac{2\,\|V\|^{(\alpha,1)}}{\alpha \, v_{\min} \, q_{\min}^{\alpha}}.
\label{finer:asym:momentum:est}
\eeqn
Lemma \ref{lem:integral:estimate} was applied in the last step.
By the same estimate, which is locally uniform in $x_0$, 
\[\ov{v}^+(x_0) = \cM^{-1} \ov{p}^+(x_0) = \cM^{-1} \lim_{t\rightarrow\infty}\tilde{p}(t)\] 
exists and is continuous in $x_0\in F^+$.\\
$\bullet$
For multi-index $\delta\in\bN_0^{2dn}$ of norm $1\le |\delta|\le k-1$ the momentum estimate 
\beq
 \| \cM^{-1} \pa^{\delta}_{x_0}(\tilde{p}(t)-\tilde{p}(0))\| \le
C_2 \|V\|^{(\alpha,k)} 
v_{\min}^{-1-|\beta|}  q_{\min}^{-\alpha-|\gamma|}
\le C_{\delta} \, v_{\min}^{+1-|\beta|}  q_{\min}^{-|\gamma|}
\Leq{finer:momentum:deriv}
is derived like the position estimate in and after \eqref{so:wie}. 
We infer that at $x_0\in F^+_{\rm loc}$ asymptotic velocity
$\ov{v}^+$ is $k-1$ times continuously differentiable.
This proves item~\ref{Statement2}.\\
As the flow $\Phi\in C^{k-1}(\bR\times P,P)$ and by Property 3.\ of Thm.\
\ref{thm:final:free}, the same statement is true for $x_0\in F^+$.
This proves item~\ref{Statement1}.
\hfill$\Box$\\[2mm]
In \cite[Lemma II.2]{He}, {\sc Herbst} noted for $n=2$ that for long
range potentials the limit
$\lim_{t\to \infty}\big(q_2(t,x)-q_1(t,x)\big)$ exists, if the
asymptotic velocities coincide. His -- perhaps astonishing -- result
immediately generalizes to the $n$--body case. To see this, we modify
\eqref{def:free}, setting
\beq \widehat{F}^\pm:= \big\{x\in \widehat{P}\mid v^{\pm} (x) \mbox {
  exists, and } v^\pm(x) \notin \Delta \big\}.
\Leq{P:hat:free:pm}

\begin{lemma}\quad \label{lem:q:minus:q}\\
  For a long range potential $V$, consider initial conditions 
  $x^{(i)}_0\equiv\big(p^{(i)}_0,q^{(i)}_0\big)\in\widehat{F}^\pm$ \
  $(i=1,2)$, whose asymptotic momenta $\ov{p}^\pm\big(x^{(i)}_0\big)$
  coincide. 
  Then
  \beq
  a^\pm:=\lim_{t\to\pm\infty} \big(q(t,x^{(2)}_0)-q(t,x^{(1)}_0)\big)
  \Leq{a:pm} 
  exists. More precisely, although the estimate 
  $p\big(t,x_0\big) - \ov{p}^\pm\big(x_0\big) = \cO(|t|^{-\alpha})$
  is in general optimal in the $t\to\pm\infty$ limit,
  \beq
  p\big(t,x^{(2)}_0\big) - p\big(t,x^{(1)}_0\big) =
  \cO(|t|^{-1-\alpha})\mbox{ and }  
  q\big(t,x^{(2)}_0\big)-q\big(t,x^{(1)}_0\big)=a^\pm+\cO(|t|^{-\alpha}).
  \Leq{pp:qq}
  Finally, if $a^\pm=0$, then $x^{(1)}_0=x^{(2)}_0$. 
\end{lemma}
\textbf{Proof:}\\
$\bullet$ To begin with, the estimate
$p(t,x_0) - \ov{p}^\pm(x_0) = \cO(|t|^{-\alpha})$
follows from \eqref{finer:asym:momentum:est} and the propagation
estimate \eqref{propagation:est}, and its optimality from
\[ \eh \big\langle \ov{p}^\pm(x_0), {\cal M}^{-1}\ov{p}^\pm(x_0)\big\rangle 
= \eh \big\langle p(t,x_0), {\cal
  M}^{-1}p(t,x_0)\big\rangle + V\big(q(t,x_0)\big). \] 
$\bullet$
The second estimate in \eqref{pp:qq} and \eqref{a:pm} follow by integration
from the first estimate in \eqref{pp:qq}.\\
$\bullet$
To derive it and the last statement, we argue like in \cite[Lemma II.2]{He}.
\hfill $\Box$\\[2mm]
We can also apply Theorem \ref{thm:both:moeller}, which is formulated for
non-singular potentials, to the unregularized Hamiltonian flow with an 
$(\alpha,k)$--potential $V:\widehat M \to \bR$.

The point is, some $x$'s end in collision, so have no well-defined
asymptotic velocity.  As in definition~\ref{def:free}, the phase space
regions $\widehat{F}^\pm\subseteq \widehat{P}$ are open. The escape
time $T^\pm (x_0)$ for initial conditions $x_0$ lying in
$\widehat{F}^\pm$ are $\pm \infty$
whereas  $T^-$ ($T^+$) are still upper (respectively lower)
semicontinuous. 

\begin{corollary}[Asymptotic velocities for singular
  potentials]\quad\label{cor:hat:free}\\[-6mm] 
\begin{enumerate}[1.]
\item 
For $\alpha>0$ and $(\alpha,k)$--potentials $V\in C^k(\widehat M,\bR)$, 
see \eqref{V:k:norm}, the restricted asymptotic velocity maps
$\ov{v}^\pm$ are $C^{k-1}$ over $\widehat{F}^\pm$. 
\item 
So for $(-\alpha)$--homogeneous potentials, that is 
\beq
\textstyle
V(q) := \sum_{1\le i<j\le n} \frac{I_{i,j}}{\|q_i-q_j\|^\alpha}.
\Leq{hom:pot}
the asymptotic velocities $\ov{v}^\pm$ are smooth on $\widehat{F}^\pm$. 
\end{enumerate}
\end{corollary}
\textbf{Proof:}
\begin{enumerate}[1.]
\item 
The flow  
is   $C^{k-1}$ on its domain and if  $x_0\in F^+$ 
then there is a time $t\ge 0$
so that $\Phi_t (x_0)\in F^+_{\rm loc}$.  (The norm $\|V\|^{(\alpha{,2})}$
appearing in Thm.\ \ref{thm:final:free} needs to be appropriately
re-defined to account for blow-up along $\Delta$).  
Then by Theorem \ref{thm:both:moeller}.1 the 
restriction of the asymptotic velocity $v^+$
to $F^+$ is a  $C^{k-1}$ with values in $\bR^{dn}$.
\item 
$V$ in \eqref{hom:pot} has finite $\|V\|^{(\alpha,k)}$ norm for any $k\in \bN$.
\hfill $\Box$%\\[2mm]
\end{enumerate}

\section{The M\o ller semi-conjugacy (short range)}
\label{sect:moeller}

We will now show that for a short range potential the flow and the
free flow are semi-conjugate, using the so-called M\o ller
transformation. 

If the potential is short range ($\alpha>1$ in \eqref{V:decay}), 
then we  can establish the asymptotics
$$q(t) = at + b + \cO\big(t^{1-\alpha}\big) \qquad \mbox{as} \ t \to + \infty$$
for forward free solutions, see equation \eqref{Moe:E} below.  
The vector $a$ is $v^+ (x_0)$ if $x_0 := (\cM\dot{q}(0), q(0))$ is the initial
condition for $q(t)$.  The vector $b$ is something like the ``impact
parameter'' found in standard treatments of classical scattering.  We
would like to think of $a, b \in \bR^{dn}$ as initial conditions at
$t = + \infty$.

One way to formalize this idea is via the M\o ller transformation,
which compares the given flow to that of a free particle.

\begin{definition}
  \label{def:moller}
  The \emph{(forward) M\o{}ller transformation}, where the (pointwise) limit
  exists, is the map
  $\Omega= \Omega_+ := \lim_{t \to + \infty} \Phi_{-t}
  \circ\Phi^0_{t}: P \dashedrightarrow P$.
  Similarly the \emph{backward M\o ller transformation} is
  $\Omega_- := \lim_{t \to - \infty} \Phi_{-t} \circ\Phi^0_{t}$, where
  the limit exists.
\end{definition} 

\begin{figure}[h]
  \begin{center}
    \includegraphics[width=100mm]{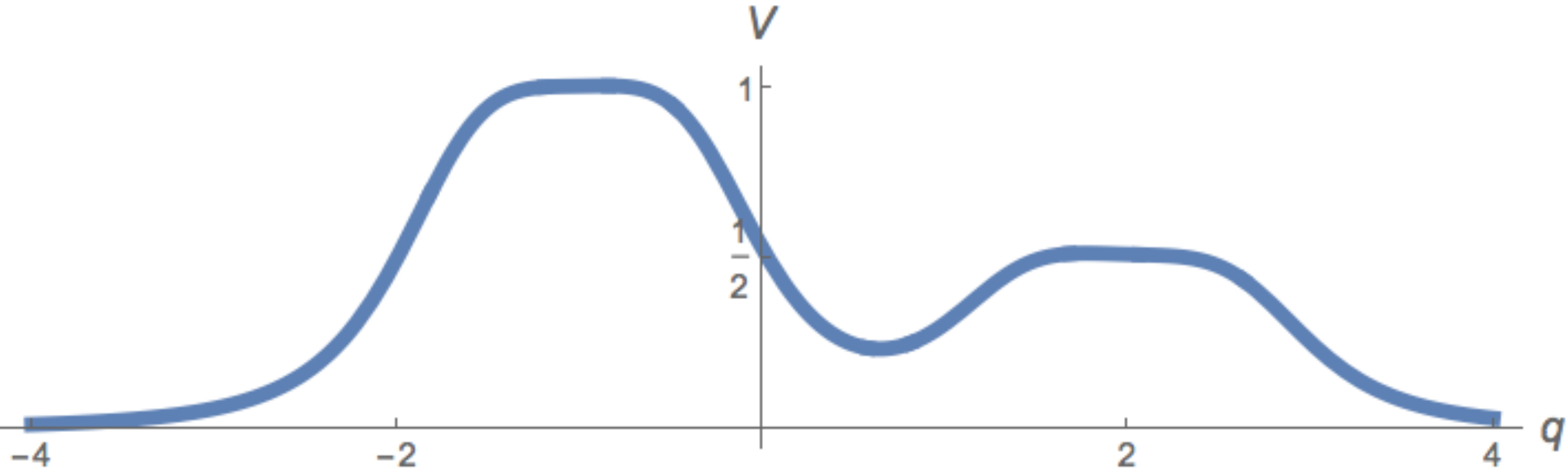}
    \includegraphics[width=100mm]{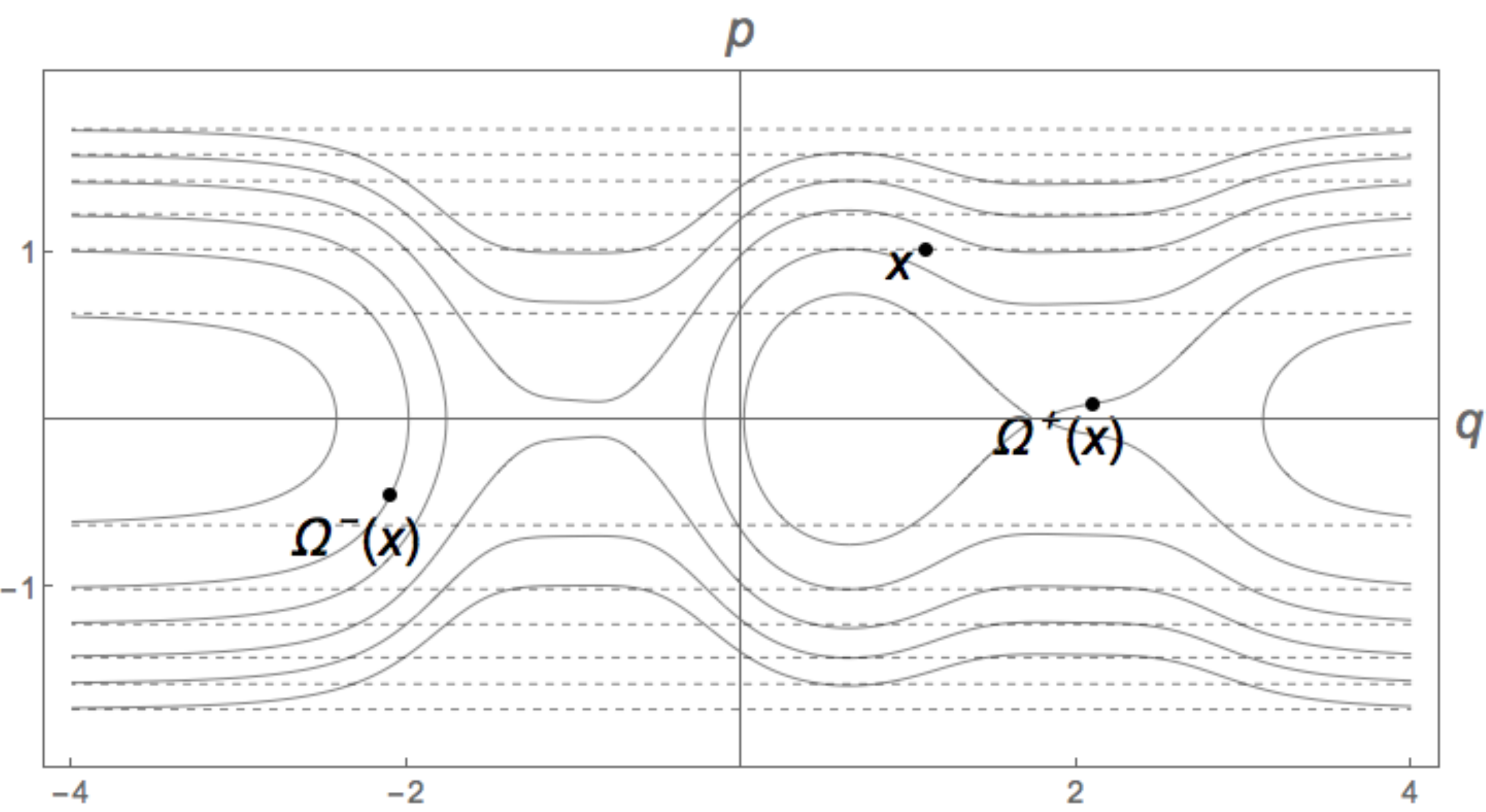}
  \end{center}
  \caption{Above: Potential $V$.  Below: Level lines of the
    Hamiltonian $H$ and of the kinetic energy $T$ (dashed); and the
    corresponding M\o ller transformations.}
  \label{fig:Moeller}
\end{figure}

See figure \ref{fig:Moeller} for a depiction of the forward and
backward M\o{}ller transformations. We have continued to use the
broken arrow notation in the definition of the M\o{}ller
transformation to allow ourselves vagueness about its domain.  We
repair this vagueness now. Moreover, these transformations provide us
with semi-conjugacy in the short range case.

\begin{theorem}[M\o{}ller transformation]\label{thm:moeller}\quad\\%
If the $(\alpha,k)$--potential $\,V$ is short range 
($\alpha>1$ in definition \ref{def:alpha-k}), then
  \begin{enumerate}[1.]
  \item \label{Statement3} 
  For $F^+_0$ and $F^+$ defined in \eqref{P0free} respectively in 
  \eqref{P:free:pm:B}, the M\o ller transformation 
    \beq \Omega = \lim_{t \to + \infty} \Phi_{-t}
  \circ\Phi^0_{t}: F^+_0 \to F^+
    \Leq{M:T}
    exists and is a $C^{k-1}$ symplectomorphism intertwining $\Phi_t$
    with $\Phi^{(0)}_t$:
    \beq 
    \Omega \circ \Phi^{(0)}_t =\Phi_t\circ\Omega \qquad (t\in \bR).
    \Leq{intertwining}
  \item  \label{Statement4}
    If $|\delta| \le k-1$,  $x_0 = (p_0, q_0) \in F^+_{\rm loc} $ 
    and $\Omega (x_0) = X_0=(P_0,Q_0)$ then the inverse M\o ller transformation $\Omega^{-1}$
    satisfies the regularity estimates: 
    \begin{align}
      \pa^\delta_{X_0} (P_0-p_0) &= 
                                   \cO\l( \|V\|^{(\alpha,k)}
                                   v_{\min}(X_0)^{-1-|\beta|}  \LA
                                   q_{\min}(X_0)\RA^{-\alpha-|\gamma|}
                                   \ri) , 
                                   \label{Moe:E:p} \\
      \pa^\delta_{X_0} (Q_0-q_0) &=
                                   \cO\l( \|V\|^{(\alpha,k)}
                                   v_{\min}(X_0)^{-2-|\beta|}  \LA
                                   q_{\min}(X_0)\RA^{1-\alpha-|\gamma|}
                                   \ri). 
                                   \label{Moe:E}
    \end{align}
  \end{enumerate}
\end{theorem}

We will deal with long range $(\alpha,k)$--potentials with $\alpha\in
(1/2,1]$ later; we will show existence of modified 'Dollard' M\o ller
transformations in Theorem \ref{thm:Dollard:Moeller}, indicating how
to generalize this to $\alpha\in (0,1]$.

\medskip
    
Let us pause to see how Theorems \ref{thm:final:free} and
\ref{thm:moeller} are related.  Suppose that
$\Omega( A, B) = x_0 \in F^+$.  We claim that $A = \cM v^+ (x_0)$
where $v^+ (x_0) \notin \Delta$ is $x_0$'s asymptotic velocity of
definition~\ref{def:asymptotic:velocity} and described in
Theorem~\ref{thm:main1}.  Inverting, we have
$\Omega^{-1} (x_0) = (A,B)$ and
$\Omega^{-1} = \lim_{t \to \infty} \Phi ^0 _{-t} \circ\Phi_{t}$.
Write $\Phi_t (x_0) = (p(t, x_0), q(t,x_0))$ and set
$p(\infty) =\cM v^+ (x_0) = A$.  Theorem \ref{thm:main1} tells us that
for $t$ large we have $p(t, x_0) = A + o(1)$.  Now the momentum $p$ is
constant under the free flow so that for large $t$ we have
$\Phi ^0 _{-t} \circ\Phi_{t} (x_0) = ( A + o(1), Q(x_0; t))$.  Letting
$t \to \infty$ kills the $o(1)$ term and yields the claim:
$\Omega ^{-1} (x_0) = (A, *)$.\\[2mm]
\textbf{Proof of Theorem \ref{thm:moeller}:}\\
$\bullet$
We now prove for $(\alpha,k)$--potentials $V$ of short range ($\alpha>1$)
pointwise existence and smoothness properties of the M\o ller transformation.\\
Thus let  $X_0=(P_0,Q_0)\in F_0$ and
 write $\tilde{Q}(t):= Q_0+\cM^{-1}P_0t$ for the corresponding free
 solution. (See \eqref{free:solu}.) Define the map $\cF_{X_0,T}\equiv$
\beq
 \cF:\widehat{\cD}  \to C\big([T,\infty), \bR^{dn}\big)\ \mbox{ , } \
(\cF r)(t) = -\cM^{-1}\!\! \int_t^\infty \!\!\! \int _s^\infty \!\!\!\!
\nabla V\big((\tilde{Q}+r)(\tau)\big)\, d\tau\, ds 
\Leq{contraction map}
on the complete metric space 
\beq
\hspace*{-2mm}
\widehat{\cD} \equiv \widehat{\cD}_{X_0,T}  := 
\Big\{ r\in C\big([T,\infty), \bR^{dn}\big)\ \Big|\
\|r\|:=\sup_{t\ge T}\|r(t)\| \le \eh q_{\min} (X_0)\Big\} .
\Leq{def:DC}
By the short range 
assumption on $V$ the  map $\cF$ is well-defined, and any function $u = \cF (r)$ in its
image satisfies  $\lim_{t\to\infty} u (t) = 0$.\\
We search for solutions $r$ of the fixed point problem $r=\cF_{X_0}(r)$.
Out of such a fixed point $r$ we will build $\Omega (X_0)$. 
First, observe that if $r$ is such a fixed point then  
\begin{equation} 
q:= \tilde{Q} + r
\label{perturbation}
\end{equation} satisfies Newton's equations 
$\ddot{q}(t)= - \cM^{-1}\nabla V\big(q(t)\big)$
and is asymptotic to $\tilde Q$.
When $X_0\in F^+_{\rm loc} \subseteq F_0$, then 
by \eqref{propagation:est} and \eqref{def:DC} the interparticle distances 
$\tilde {q}_{i,j}(\tau)\ge \eh(q_{i,j}+ v_{i,j}\tau)$. Thus, using
\eqref{V:k:norm} 
\[ \big\| \cM^{-1}\nabla V\big((\tilde{Q}+r)(\tau)\big) \big\| \
\le \ \|V\|^{(\alpha)}\, \LA \eh(q_{\min} + v_{\min}\tau)\RA^{-1-\alpha}
\qquad (\tau\ge0).\]
So by Lemma \ref{lem:integral:estimate}
$\big\|  \cM^{-1}\int_s^\infty \nabla V\big((\tilde{Q}+r)(\tau)\big) \,d\tau \big\| 
\le \frac{2 \|V\|^{(\alpha)}}{v_{\min}}  \LA\eh(q_{\min} + v_{\min} s)\RA^{-\alpha}$
and 
$\|(\cF r)(t)\| \le \frac{2\, \|V\|^{(\alpha)}}{(\alpha-1)v_{\min}^2} 
 \LA\eh q_{\min}\RA^{1-\alpha}
\le  \frac{8d\|V\|^{(\alpha,1)}}{(\alpha-1)v_{\min}^2}  (\eh q_{\min})^{1-\alpha}
\le \eh q_{\min}$,\\ 
as $X_0 \in F^+_{\rm loc}$ and $\delta\le \alpha-1$.
So $\cF$ maps $\widehat{\cD}$ into itself.

Next we show that $\cF$ is a contraction on $\widehat{\cD}$. 
So let $r^{(0)}\neq r^{(1)}\in \widehat{\cD}$. Then 
\[ \frac {\big\|\cF\big(r^{(0)}\big)-\cF \big(r^{(1)}\big)\big\|} { \|r^{(0)}-r^{(1)}\| }
\le \int_0^\infty \!\!\! \int _s^\infty \!\!\!
\int_0^1\big\| \cM^{-1} D\nabla V\big((\tilde{Q}+r^{(\rho)})(\tau)\big) \big\|
\,d\rho \; d\tau\; ds\]
with $r^{(\rho)}:=(1-\rho)\,r^{(0)}+\rho\, r^{(1)}$.
The right hand side is majorized by
\[  \|V\|^{(\alpha,2)}\int_t^\infty \!\!\! \int _s^\infty 
\LA\eh(q_{\min}+ v_{i,j} s)\RA^{-2-\alpha} \,d\tau\; ds
\le  \frac{2\|V\|^{(\alpha,2)}}{\alpha(1+\alpha) v_{i,j}^2 q_{\min}^\alpha}
\le \frac{\delta}{16dn} < 1. \]
By Banach's theorem $\cF_{X_0}$ has a unique fixed point $r$. Evaluating
the corresponding solution $\tilde Q (t) + r(t)$ to Newton's equations appropriately
at $t = 0$   yields
the value of the M\o ller transformation on $X_0$.
Indeed we claim that  
\[  \Omega (X_0) = \big(P_0+\cM \dot{r}(0),\, Q_0+r(0)\big) . \] 
To see this, we approach the problem of approximating $r(t)$ ``from the other end of time''
as follows.  Write $\Phi_{-T}\circ\Phi^{(0)}_T(X_0) = \big(\cM (\dot r^{(T)})(0),\, r^{(T)}(0) \big)$.\\
Then the solution $r^{(T)} :[0,T] \to \bR^{dn}$  to Newton's equations
with  initial position  $r^{(T)}(0)$ and initial velocity $\dot r^{(T)}(0)$ 
is the unique fixed point of the map 
\[\cF^{(T)}:\widehat{\cD}^{(T)}  \to \widehat{\cD}^{(T)} 
\mbox{ , }\,
(\cF^{(T)} r^{(T)})(t) = -\cM^{-1}\! \int_t^T \!\!\! \int _s^T \!\!\!
\nabla V\big((\tilde{Q}+r)(\tau^{(T)})\big) \,d\tau\; ds\]
on 
\[\widehat{\cD}^{(T)}  := \l\{ r\in C\big([0,T], \bR^{dn}\big) \mid 
\max_{t\in[0,T]}\|r(t)\| \le \eh q_{\min} (X_0), r(T)=\dot{r}(T)=0\ri\} \!,\]
and by uniqueness of the original fixed point $r$  we must  have that 
\[r(t)=\lim_{T\to+\infty} r^{(T)}(t)\qmbox{,}\dot{r}(t) =
\lim_{T\to+\infty} \dot{r}^{(T)}(t) \qquad (t\ge0).\]
To see that M\o ller transformation is defined on all of $F_0$,
observe that for any $X_0\in F^+$ we have, eventually, for large
enough times $h$ that $\Phi^{(0)}_h(X_0)\in F_{\rm loc}$, at which point
we have just seen that $\Omega\big(\Phi^{(0)}_h(X_0) \big)$ exists.  Then
observe by inspecting the definition of the limits that
$\Omega (X_0) = \Phi_{-h} \circ \Omega \circ \Phi^{(0)}_h (X_0)$. \\ 
As a locally uniform limit the M\o ller transformation is continuous on
$F^+$.  The intertwining relation \eqref{intertwining}
follows, since the flows are $\bR$--actions, or alternatively by
re-arranging the just-proved relationship,
$\Omega = \Phi_{-h} \circ \Omega \circ \Phi^{(0)}_h $ valid for all
sufficiently large $h$ in a neighborhood of any $X_0$.\\[2mm]
$\bullet$
To investigate the degree of smoothness of $\Omega^+$,
instead of the operator \eqref{Q:op} related to the initial value problem, 
we now use the operator $\cP\equiv \cP_{X_0}$, with
\beq
\cP(w)(t) := 
-\cM^{-1}  \int_t^\infty \!\! \int_s^\infty  D\nabla V\big(\tilde{Q}(\tau)\big)\, w(\tau)\,d\tau \,ds
\qquad(t\geq 0),
\Leq{P:op}
on the Banach space $C^b([0,\infty),\bR^{dn})$ of bounded curves.
Its operator norm is majorized by 
\begin{align}
\|\cP_{X_0} \| 
&\le\|V\|^{(\alpha,2)}\int_0^\infty \!\! \int_s^\infty \!\! \LA 
\tilde{Q}(\tau)\RA^{-2-\alpha} d\tau \,ds \NN\\
&\le \|V\|^{(\alpha,2)}\int_0^\infty \!\! \int_s^\infty \!\! \
\LA \eh(q_{\min} + v_{\min} s) \RA^{-2-\alpha}\! d\tau \,ds \NN\\
&\le \frac{2^{2+\alpha}}{\alpha} \|V\|^{(\alpha,2)} q_{\min}^{-\alpha} v_{\min}^{-2}
\le \frac{2^{2+\alpha}}{\alpha}\frac{\delta}{16dn} < 1 \NN
\end{align}
if $\alpha \le3$ (for larger $\alpha$ one uses the forward flow into
$F^+_{\rm loc}$, where the estimates become better). So we can invert
${\rm Id} - \cP_{X_0}$ in order to solve for $|\delta| \le k-1$ 
\beqn
\lefteqn{\pa^\delta_{X_0} r(t,X_0) = -\int_t^\infty \!\! \int_s^\infty 
 \cM^{-1} \pa^\delta_{X_0}\nabla V\big(q(\tau,X_0) \big)\,d\tau \,ds 
= - \! \sum_{N=1}^{|\delta|}  \cM^{-1} \, \times} &&
\label{r:est}\\
&&
\hspace*{-6mm}
\times\hspace*{-4mm} \sum_{\stackrel{\delta^{(1)}+\ldots+\delta^{(N)}
  = \delta} {|\delta^{(i)}|>0}} 
 \int_t^\infty \!\! \int_s^\infty\! D^N\nabla V\big(q(\tau,X_0)\big)
\l(\pa^{\delta^{(1)}}_{X_0}q(\tau,X_0),\ldots, 
\pa^{\delta^{(N)}}_{X_0}q(\tau,X_0)\ri) d\tau \,ds \NN
\eeqn
with the shorthand $q=\tilde{Q}+r$ in a way similar to \eqref{oneminq}.
This shows \eqref{Moe:E} and finishes the proof of Claim \ref{Statement4}.\\ 
As $C^1$--limit of the symplectomorphisms
$\Phi_{-t}\circ\Phi^{(0)}_t$ the M\o ller transformation $\Omega^+$ is a
symplectomorphism onto its image.
But this image coincides with $F^+$, by its mere definition \eqref{P:free:pm:B} and by reversing the
roles of the two flows.\\
So Claim \ref{Statement3} is also true.
\hfill $\Box$
\begin{remark}[M\o ller transform]\quad\\
  The standard reference for the {\em M\o ller} transform is section 5
  of \cite{DG} by {\sc Dere\-zi\'{n}ski} and {\sc G{\'e}rard}.
  In the case of finite-range interactions {\sc Hunziker}, in
  \cite{Hu1,Hu2} proved that the M\o ller transform exists and used it
  to establish {\em asymptotic completeness} of finite range
  interactions.  This asymptotic completeness includes the
  decomposition of solutions into independent clusters where `cluster'
  has the meaning alluded to above.  
  
  Hunziker viewed the M\o ller
  transform as the classical version of the quantum M\o ller
  transform, or wave map, defined as the limit of
  $\exp(-it H) \exp(it H_0)$ as $t \to \infty$. Here $H = H_0 + V$ is
  the quantum version of our Hamiltonian so that $H_0$ corresponds to
  a multiple of the Laplacian on $\bR^{dn}$.\\ 
  Soon afterwards, {\sc
    Simon} \cite{Sim} used the method to establish asymptotic
  completeness for the classical two-body problem with short range
  interactions provided the second derivative of the potential decays
  appropriately. In an appendix Simon exhibited the necessity of his
  second derivative decay conditions by constructing a potential for
  which his decay conditions failed and which admits two distinct
  hyperbolic solutions asymptotic to the same free solution. Thus
  $\Omega^{-1} (x_0) = \Omega^{-1} (y_0)$ for $x_0, y_0$ not lying on
  the same orbit, so that whatever $\Omega$ is, it is at least
  ``two-valued'' and not a well-defined map.\\ 
  {\sc Derezi\'{n}ski} and {\sc G{\'e}rard}, among many other results, established the existence
  and invertibility of the M\o{}ller transformation for potentials of
  superexponential decrease in \cite[sect.\ 5.10]{DG}. \hfill $\Diamond$
\end{remark}

% ----------------------------------------------------
\section{The Dollard-M\o ller semi-conjugacy (long range)}
\label{sect:moeller2} 
%----------------------------------------------------

The gravitational and Coulomb  potentials are long range but not short
range so the M\o ller transformation fails to exist for them.  Dollard \cite{Do} discovered that  by modifying the comparison free dynamics
 in a time-dependent way he could define a modified   M\o
ller transformation which  existed for  long range potentials.  
We will call his modified transformation  the Dollard-M\o ller
transformation.   It  will yield the asymptotics 
\beqno
p(t)  = \cM v + o(1)\\
 q(t)  = v t + W(t, \cM v) + b + o(1) \\
 \text{ with }  v= v^+ (x(0)) \mbox{ and }   W(t, \cM v) = o(t)  \mbox{ as } \ t \to + \infty
 \label{eq:Chazylike}
\eeqno
valid for all escape solutions  $x(t) = (p(t), q(t))$ and all   long-range potentials ($0 < \alpha \le 1$ in \eqref{V:decay}),
whether they have singularites or not. 
See equation \eqref{def:W} for the relation between $W$ and the potential $V$. 
This assertion on asymptotics  follows from the existence of the inverse Dollard-M\o ller transformation $\Omega^{-1}$,   part 1 of \ref{thm:Dollard:Moeller}.
See remark \ref{rmk:asymptotics}  for a sketch of a proof of a derivation of \eqref{eq:Chazylike} from part 1. 
The asymptotic velocity $v$ occurring in    the asymptotics \eqref{eq:Chazylike} is given by   $\Omega^{-1} (x(0)) = (v, \beta)$ for some $\beta$.
The ``impact parameter'' $b$, projected onto $v^{\perp}$ represents the affine orbital parameter
described in  part 3 of \ref{thm:Dollard:Moeller} below.

\begin{definition}\label{def:dollard}\quad\\
The Dollard dynamics $\Phi ^D _{t, s}$ (see~\eqref{groupoid}) associated with a  potential $V$ on $\bR^{dn}$
  is the non-autonomous flow defined by the time dependent
  \emph{Dollard Hamiltonian} 
  \[ H^D := K+\tilde{H}^D:\ \bR_t\times F_0
  \longrightarrow \bR   \]
  given by 
  \beq
  H^D (t, p, q) = \eh \langle p, \cM^{-1}p \rangle + V \big( \langle
  t\rangle\cM^{-1}p \big)  \text{ where}  \LA t\RA=\sqrt{1+t^2} 
    \Leq{dollard:hamiltonian2}
      \end{definition}
      
The first term  $K$ of $H^D$ is the usual kinetic energy. Its second term 
  $\tilde{H}^D_t(p,q)$ is the potential turned into a function of momentum.
$H^D$ is independent of $q$ so the  momentum $p$ is constant along the non-autonomous Dollard
flow $\Phi ^D _{t, s}$.

\begin{example}[Newtonian case] \quad\\
  Take the case of the Newtonian $n$-body problem, where the potential
  is homogeneous of degree $-1$.  Using
  $\langle t\rangle  = t( 1 + \frac{1}{2} \frac{1}{t^2} + \ldots)$ for $t \gg 1$ we
  see that
  $H_D = \frac{1}{2} \langle p, \cM^{-1} p\rangle + V( \langle t
  \rangle \cM^{-1} p) = \frac{1}{2} \langle p, \cM^{-1} p\rangle +
  \frac{1}{t} V( \cM^{-1} p) + \cO(1/t^3)$
  for large $t$, where the $\cO(1/t^3)$ term depends only on $p$. Then
  the ODEs to solve to find the Dollard flow are
  $$
  \begin{cases}
    \dot q = \cM^{-1} p + \frac{1}{t} \nabla V( \cM^{-1} p) + \cO(1/t^3)\\
    \dot p = 0
  \end{cases}
  $$
  which integrate to yield precisely Chazy's asymptotics
  ~\eqref{hyperbolicChaz}  above.   Compare with 
  {\sc Chazy} \cite[page 46]{Cha},
  1922. One could argue that the proper Dollard Hamiltonian
  \eqref{dollard:hamiltonian:homogeneous} has Chazy's work as a  precursor. 
  \hfill $\Diamond$
\end{example} 

Returning to a general $V$, we compute the time-dependent flow of $H^D$, for
initial time $s\in\bR$ and final time $t\in\bR$, to have the form: 
\beq
\Phi^D_{t,s}(p,q) = \big(\,p\; ,\; q+(tv+ W(t;p))-(sv+ W(s;p))
\big) \quad \big( (p,q)\in F_0 \big) .
\Leq{Dollard:dynamics} 
where  \beq W: \bR_t\times F_0\to
\bR^{dn}\qmbox{,} W(t;p) = \int_0^t \nabla_p V\big(\langle s\rangle
\cM^{-1}p\big)\,ds\,. 
\Leq{def:W} 
If $V$ is an $(\alpha, k)$ potential then 
$W\in C^{k-1}\big(\bR_t\times F_0,\bR^{dn}\big)$,  
and
\[t\mapsto W(t;p) =\l\{
\begin{array}{ll}
\cO(|t|^{1-\alpha})&,\, \alpha\in (1/2,1)\\[1mm]
\cO(\log(|t|))&,\, \alpha=1,
\end{array}
\ri.  \qquad(|t|\to\infty).\]
Although the correction term $W(t,p)$ to linear motion can go to infinity with $t$,
we have that $W(t, p) = o(t)$, which is to say, that  $|t|  > > |W(t;p)|$ as $t \to \infty$. 
It will be crucial below that  for fixed $p$ and $t_0$,  
\beq W(t + t_0, p) - W(t, p) \to 0 \text{ as } t \to \infty \, ,
\Leq{W:differences} 
as the reader can easily verify. 

The asymptotic velocity of any Dollard solution curve $\Phi^D_{t,s}(x_0)$ with  
$x_0 = (p,q)$ is $v=\cM^{-1}p$. All Dollard solutions \eqref{Dollard:dynamics} 
which share a fixed initial momentum $p$ are translates of one another:
\[\Phi^D_{t,s}\big(p,q^{(2)}\big)-\Phi^D_{t,s}\big(p,q^{(1)}\big) = \big( 0,q^{(2)}-q^{(1)} \big) \qquad 
\big( s,t\in\bR,\, q^{(i)}\in\bR^{dn} \big).\]

For explicit computations of Dollard flows and comparison of the induced transformations
with M\o ller transformations  see the appendices. 

We will  use the Dollard flow $ \Phi^{D}_{0,t}$  in place of the free flow $\Phi^{(0)}_t$ in
order to define a version of the M\o ller transformation. However, collisions in backward time prevent 
us from defining a direct Dollard-M\o ller transform on $\widehat{F}^+$.  
The backward  $n$-body flow $\widehat \Phi_{-t}$,  
 $t > 0$,  applied to some   points of $\widehat{F}^+$ 
may not exist  due to multi-body collisions in backwards time.    
To circumvent this problem we instead define the   {\em inverse}  Dollard-M\o ller transformation, 
whose definition only uses the forward flow, so that its domain  can be
taken to be $\widehat{F}^+$.  

\begin{theorem}[Dollard-M\o ller transformations]\quad
  \label{thm:Dollard:Moeller}\\
  (The reader may wish to refer to subsection \ref{subsec:notations} for notations.) 
  For long range $(\alpha,k)$--potentials $V$ (see \eqref{V:k:norm})
  with $\alpha\in (1/2,1]$  and collision singularities allowed, the following hold.  
  \begin{enumerate}[1.]
  \item  \label{long:range:1}
  The backward and forward  inverse  Dollard-M\o ller transformations 
    \beq
    \Omega^{-1, \pm} :=  \lim_{T\to\pm\infty}
    \Phi^{D}_{0,T}\circ\widehat{\Phi}_{T} 
    \qmbox{,}
    \Omega^{-1, \pm}: \widehat{F}^\pm\to F_0
    \Leq{I:M:T:D}
    exist in the sense of locally uniform convergence.
  \item \label{long:range:2}
  a) These transformations   conjugate the n-body flow on  $\widehat{F}^\pm$
  with   the  free flow 
    \beq
    \Omega^{-1, \pm}\circ \widehat{\Phi}_t =  \Phi^{(0)}_t\circ
    \Omega^{-1, \pm} \,.
    \Leq{conjugate}
    b) For $k\ge3$ the $\Omega^{-1, \pm}$ are
    $C^{k-2}$--smooth symplectomorphisms onto their images. 
  \item \label{long:range:3} The analog of \eqref{Moe:E:p}  holds for
    $\Omega^{-1, \pm} - Id$.
  \item \label{long:range:4} 
    For any   $v\in \bR^{dn}\,\backslash\,\Delta$,
    the space of orbits having  asymptotic velocity $v$ form an affine space with 
    underlying vector space the  tangent space of the sphere $S^{dn-1}$ at
    $v/\|v\|_\cM$.
  \end{enumerate}
 \end{theorem}
\textbf{Proof:}\\ % of Theorem \ref{thm:Dollard:Moeller}
We will make use of the open subset ${F}_{\rm loc}^+$ of
$\widehat{P}$ defined by precisely the same conditions as
${F}_{\rm loc}^+$ in \eqref{def:ff} with all points  lying in $\widehat{P}$ - the
phase space points with no  collisions. Note that
for $\alpha$-homogeneous potentials, the  conditions within
\eqref{def:ff}  respect the homogeneity of kinetic and potential
energy. \\
$\bullet$
As both flows $\widehat{\Phi}$ and $\Phi^{D}_{\bullet,\bullet}$ are
$C^{k-1}$--smooth on their  maximal domains $\widehat{D}$ and
$\bR_t\times\bR_s\times F_0$ respectively, by Theorem
\ref{thm:final:free}.3 and its Corollary \ref{cor:hat:free} we can
assume without loss of generality that $x_0=(p_0,q_0)\in
\widehat{F}_{\rm loc}^+$. 

We consider the Dollard solution \eqref{Dollard:dynamics}, 
$t\mapsto \Phi^{D}_{t,0}(X)$ 
with initial value $X\in F_0$ and
denote by $X_T(x_0)$ the initial value  with the property
\[\Phi^{D}_{T,0}\big(X_T(x_0)\big) = \widehat{\Phi}_{T}(x_0) \qquad 
\big( x_0\in \widehat{F}_{\rm loc}^+,\ T\ge0 \big).\]
Since $\Phi^{D}_{\bullet,\bullet}$ is the solution of a time dependent
initial value problem, we have  
\beq
\Phi^{D}_{t,t}={\rm Id}_{F_0}\qmbox{and}
\Phi^{D}_{t_2,t_1}\circ \Phi^{D}_{t_1,t_0} = \Phi^{D}_{t_2,t_0}\qquad
(t,t_i\in\bR), 
\Leq{groupoid}
so that $\big(\Phi^{D}_{t_1,t_0})^{-1}= \Phi^{D}_{t_0,t_1}$. \\

\noindent{\large {\it Proof of part  1 of the theorem.}}\\ 
In \eqref{I:M:T:D}  we claim  pointwise existence and local uniformity of the limit $T \to \infty$ of   
\beq
\Omega^{-1}_T := (\Phi^D_{T,0})^{-1}\circ\widehat{\Phi}_T =
\Phi^D_{0,T}\circ\widehat{\Phi}_T. 
\Leq{Om:star:T}
We compute, with  $v:={\cal M}^{-1}p$ denoting velocity, that 
\begin{align} 
\Omega^{-1}_T(x_0) 
&= \textstyle
\big(p(T,x_0), q(T,x_0)-v(T,x_0)T-\int_0^T \nabla_p V\big(\langle s\rangle \cM^{-1}p(T,x_0)\big)\,ds\big)\NN\\
&=
\big(p(T,x_0) , q_0+ r(T,x_0)\big), 
\label{Omega:T}
\end{align}
where 
\begin{align} 
\label{r:for:Omega}
 r(T,x_0) := &\,\textstyle
\int_0^T \big[v(s,x_0) \!-\!v(T,x_0)\! -\!\nabla_p V\big(\langle s\rangle \cM^{-1}p(T,x_0)\big)\big] \,ds \big) \\
 = &\,\textstyle
 \cM^{-1}\! 
 \int_0^T \big[\int_s^T  \nabla V\big(q(\tau,x_0)\big) d\tau -
 \!\langle s\rangle \nabla V\big(\langle s\rangle \cM^{-1}p(T,x_0)\big) \big]\,ds,\NN
 \end{align} 
see \eqref{Dollard:dynamics} and \eqref{def:W}.\\
$\bullet$ 
We begin the proof of \eqref{I:M:T:D}  by  showing that
\[r^+(x_0):=\lim_{T\to+\infty} r(T,x_0)=
r(0,x_0)+\lim_{T\to+\infty}{\textstyle \int_0^T}\dot{r}(t,x_0)\,dt \] 
exists.
Therefore, we first estimate its $T$--derivative.
\begin{align} 
\label{eq:dot:r}
\hspace*{-4mm}
\dot{r}(T,x_0) = &
                   - T\,\dot{v} (T,x_0) - \langle T\rangle \nabla
                   V\big(\langle T\rangle \cM^{-1}p(T,x_0)\big) \\  
& +\cM^{-1} \textstyle \int_0^T 
\langle s\rangle^2 D\nabla V\big(\langle s\rangle
  \cM^{-1}p(T,x_0)\big) \,ds \; \cM^{-1}\nabla V\big(q(T,x_0)\big)\,.\NN
\end{align} 
The propagation estimate \eqref{propagation:est} and \eqref{finer:asym:momentum:est} imply that
$\|p(T,x_0)-p^+(x_0)\| = \cO\big(\langle T\rangle^{-\alpha}\big)$ and thus 
locally uniformly in $x_0\in \widehat{F}_{\rm loc}^+$
\[\|q(T,x_0)-(q_0+v^+(x_0)T)\| = \left\{\begin{array}{cl}
\cO\big(\langle T\rangle^{1-\alpha} \big) &,\, \alpha < 1 \\[1mm] 
\cO\big(\log(T)\big) &,\, \alpha = 1
\end{array}
\right. .\]
\begin{enumerate}[1.]
\item 
For $\alpha\in(1/2,1)$ the first line on the right hand side of \eqref{eq:dot:r} equals
\begin{align*}
&\cM^{-1}  \big[T\,\nabla V(q(T,x_0)) - \langle T\rangle  \nabla V\big(\langle T\rangle \cM^{-1}p(T,x_0)\big) \big]\\
=&\,
T \cM^{-1}  \big[\nabla V(q(T,x_0)) -  \nabla V\big(\langle T\rangle\, v(T,x_0)\big) \big] 
+\cO\big(\langle T\rangle^{-2-\alpha}\big)\\
=&\,
T \cM^{-1}  \big[\nabla V\big(\langle T\rangle\, v(T,x_0)\! +\! \cO\big(T^{1-\alpha}\big)\big) - 
 \nabla V\big(\langle T\rangle\, v(T,x_0)\big) \big] + \cO\big(\!\langle T\rangle^{-2-\alpha}\big)\\
=&\,\cO\big(\langle T\rangle^{-2\alpha}\big) + \cO\big(\langle T\rangle^{-2-\alpha}\big)
  =\cO\big(\langle T\rangle^{-2\alpha}\big)\, ,
\end{align*}
since $\langle T\rangle-T=\cO\big(\langle T\rangle^{-1}\big)$.
\item 
As $D\nabla V\big(\langle s\rangle \cM^{-1}p(T,x_0)\big) =\cO\big(\langle s\rangle^{-2-\alpha}\big)$, 
for $\alpha\in(1/2,1)$  the second 
line of \eqref{eq:dot:r} has the order $\cO\big(\langle T\rangle^{1-\alpha}\langle T\rangle^{-1-\alpha}\big)
=\cO\big(\langle T\rangle^{-2\alpha}\big)$, too.
\item 
For $\alpha=1$ the orders of both lines in \eqref{eq:dot:r} are 
$\cO\big(\langle T\rangle^{-2}\log(\langle T\rangle)\big)$.
\end{enumerate}
We conclude that \eqref{eq:dot:r} is of order $\cO\big(T^{-2\alpha}\big)$ for $\alpha\in (1/2,1)$, respectively 
$\cO\big(T^{-2}\log(T)\big)$ for $\alpha=1$.
By our assumption $2\alpha>1$ we finally obtain existence of $r^+(x_0)$, and thus of
inverse Dollard-M\o ller transformation $\Omega^{-1} = \lim_{T\nearrow
  +\infty}\Omega^{-1}_T$.\\ 
$\bullet$
As $r^+(x_0) = \lim_{T\to+\infty} r(T,x_0)$ exists, by the analogs of \eqref{r:for:Omega} and \eqref{eq:dot:r}
\begin{align} 
 r_i&(T,x_0) = r^+_i(x_0) - \int_T^\infty \dot{r}_i(\tau , x_0)\,d\tau 
 = r^+_i(x_0) +\frac{1}{m_i} \sum_{j\in N\setminus\{i\}} \label{r:int:eq}\\
&\!\!\!\! \Big[ \int_T^\infty \!\! \Big( \tau \nabla V_{i,j}\big(q_i(\tau,x_0)-q_j(\tau,x_0)\big)
- \langle \tau\rangle
\nabla V_{i,j} \big(\langle \tau \rangle (v_i(\tau,x_0)-v_j(\tau,x_0))\big) \Big)\,d\tau  \NN\\
 &\!\!\!\!- \!\!\int_T^\infty \!\int_\tau^\infty \!\!  \langle s\rangle^2
D \nabla V_{i,j} \big(\langle {s} \rangle 
({v}_i(\tau,x_0)- {v}_j(\tau,x_0))\big)
\big(\dot{v}_i(\tau,x_0)- \dot{v}_j(\tau,x_0)\big)\, ds\,d\tau \Big].\NN
\end{align} 
When one substitutes the argument $q_i(\tau,x_0) - q_j(\tau,x_0)$ in the second line of \eqref{r:int:eq},
using $q(\tau,x_0) =  v(\tau,x_0)\tau + W(\tau;\Phi_\tau(x_0)) - r(\tau,x_0)$, then one obtains an
integral equation for $r$. 

When we assume that $r$ belongs to the complete metric space 
$\widehat{\cD}_{X_0,T}$ defined in \eqref{def:DC}, then the integrand is of order 
$\cO\big(\tau^{-2\alpha})$ for $\alpha\in (1/2,1)$ and $\cO\big(\tau^{-2}\log(\tau))$ for $\alpha=1$. 
So from \eqref{r:int:eq} we infer that  $r(T,x_0) - r^+(x_0)$ is of order $\cO\big(T^{1-2\alpha})$, resp.\ $\cO\big(T^{-1}\log(T))$.
As a function of $r$, the right hand side of \eqref{r:int:eq} is a contraction for $T$ large, 
justifying the assumption $r\in \widehat{\cD}_{X_0,T}$.
\\
$\bullet$
As convergence is locally uniform on $\widehat{F}^+$,
by the parametrized fixed point theorem the dependence of $r$ on $x_0$ is continuous.
So  the map  $r^+: \widehat {F}^+\to\bR^{dn}$ is continuous, too.

Estimates of the derivatives w.r.t.\ 
this initial condition proceed like in the proof for the short range case, 
that is, Theorem \ref{thm:both:moeller}.3.

As stated in Corollary \ref{cor:hat:free}, for $(\alpha,k)$--potentials asymptotic 
velocity $\ov{v}^+\in C^{k-1}\big(\widehat{F}^+,\bR^{dn}\big)$. 
So by \eqref{Omega:T},
$\Omega^{+,*}$ is continuous, and as smooth as $r^+$.

Note, however, that in \eqref{r:int:eq} the second derivative of the long range potential $V$ appears. 
This is different from the case \eqref{contraction map} of short range potentials, where only the
first derivative is needed.
Therefore, in comparison with Part \ref{Statement3} of Theorem \ref{thm:both:moeller}, 
we lose one derivative in Part \ref{long:range:2} of Theorem \ref{thm:Dollard:Moeller}.\\
$\bullet$
By Lemma \ref{lem:q:minus:q}, $\Omega^{-1, \pm}$ is one to one. So we can invert $\Omega^{-1, \pm}$
on its image, yielding the M\o ller transformation $\Omega^{\pm}$. We still have to prove that for any 
$x_0 = (p_0,q_0)\in \widehat{F}_{\rm loc}^+$ and its image $X\equiv X(x_0)
:= \Omega^{+,*}(x_0)$ the M\o ller transformation is of the form
\beq
\Omega^+(X) = \lim_{T\to +\infty}\Omega_T(X)\qmbox{ for}
\Omega_T := \widehat{\Phi}_{-T} \circ\Phi^{D}_{T,0}.
\Leq{Omega_T}
But this means to control $r$ as a function of $X$ instead of $x_0$.
So the analysis is similar, and we omit it.\\
{\it This completes the proof of item (1) of the proposition, i.e of~\eqref{I:M:T:D}.
}

\vskip .2cm 

\noindent{\large {\it Proof of part 2 of the theorem.}}\\ 
The intertwining property \eqref{conjugate} follows by first noting that for $\Omega^{-1}_T$ from \eqref{Om:star:T}
\[\Omega^{-1}_T \circ \widehat{\Phi}_t = (\Phi^D_{0,T}\circ \Phi^D_{T+t,0})\circ \Omega^{-1}_{T+t} \]
follows by applying the groupoid property \eqref{groupoid}, and by \eqref{Dollard:dynamics}, 
\[ \Phi^D_{0,T}\circ \Phi^D_{T+t,0}(p,q) = \big(p,q+tv +W(T+t;p,q)-W(T;p,q)\big).\]
Then $\lim_{T\to+\infty} \Phi^D_{0,T}\circ \Phi^D_{T+t,0} = \Phi^{(0)}_t$, since using \eqref{def:W}
\[\lim_{T\to+\infty}\big( W(T+t;p,q)-W(T;p,q) \big) 
= \lim_{T\to+\infty}\int_{T}^{T+t} \nabla_p V\big(\langle s\rangle \cM^{-1}p\big)\,ds =0.
\]
$\bullet$
As a locally uniform limit of symplectomorphisms $\Omega_T$ in $C^{1}$ norm, 
for $k\ge3$ the Dollard-M\o ller transformation $\Omega^+$ is a
symplectomorphism onto its image. This is shown by suitably modifying
the proof of Theorem \ref{thm:both:moeller}.3. \\
{\it This completes the proof of item 2 of the proposition.}

\vskip .2cm 

\noindent{\large {\it Proof of part 3 of the theorem.} }\\
The analog of \eqref{Moe:E:p} follows from \eqref{O:p}, as 
the Dollard dynamics \eqref{Dollard:dynamics} conserves momentum.
{\it This completes  the proof of item 3.}\\ 

\noindent{\large {\it Proof of part 4 of the theorem.} } \\
The proof relies on Lemma \ref{lem:orbits:with equal:as:ve} below,  
the conjugacy relation  \eqref{conjugate} which forms part 3 just proved, 
 and the relation \eqref{eq:meaningOmega}
proved below. 

Let us write ${\cP}_v$ for the space of all trajectories $x(t)$ having $v^{+} (x(t))= v$
where $v \notin \Delta$ is fixed.   Let $\pi^{\perp}: \bR^{dn} \to v^{\perp}$ be the orthogonal projection
so that $\pi^{\perp} (w) = w - (v \langle w, v \rangle  / |v|^2)$.  Define a map
$${\cP}_v \times {\cP}_v \to v^{\perp}$$
\beq
 (x, x^{(0)}) \mapsto   \lim_{t \to \infty}  \pi ^{\perp} ((q(t) - q^{(0)} (t)) =: b(x,x^{(0)}) \in v_* ^{\perp} 
 \Leq{AFFINE}
 where we've written by
$x(t) = (p(t),q(t)), x^{(0)} (t) = (p^{(0)}(t), q^{(0)}(t))$ for two trajectories, i.e.   points in ${\cP}_v$.   
By Lemma \ref{lem:orbits:with equal:as:ve} this limit exists and is independent of
where we start on the orbits: shifting $x(t)$ to $x(t + t_1)$ and $x^{0}(t)$ to $x^{(0)} (t + t_0)$
yields $ \lim_{t \to \infty}  (q(t+  t_1) - q^{(0)} (t+t_0) )= \lim_{t \to \infty}  (q^{(1)}(t) - q^{(0)} (t) ) + (t_1 - t_0)  v $
so leaves  the map \eqref{AFFINE} unchanged.

Think of one of the orbits, $x^{(0)}$, as the ``origin'' of ${\cP}_v$.  Then {\em we must show  that the map
\eqref{AFFINE},  viewed as a function of  $x$ alone, is onto, and that its image uniquely determines $x$
up to a time translation}.  

It will be important to understand that $x \in {\cP}_v$ iff $\Omega^{-1}(x(0)) =(\cM v, \beta)$ for some $\beta$.
This is an immediate consequence of 
\beq 
v^+ (x(t)) = v_* \iff \Omega^{-1} (x(0)) = (\cM v_* , \beta), \text{ some } \beta
\Leq{eq:meaningOmega} 
valid for all escape orbits $x(t)$. 
To establish the validity of \eqref{eq:meaningOmega} recall that the free flow
(or the Dollard flow)  does not change the momentum component. 
Write $\Omega^{-1} (x(0)) = (\cM v, \beta)$, for some $v, \beta$.  
Let ${\rm pr}_1$ denotes the projection onto the momentum factor. Then we have   
\[\cM v =  {\rm pr}_1 \Omega^{-1}(x(0)) = {\rm pr}_1  \Omega^{-1} (x(t)) = 
\lim_{t \to \infty} {\rm pr}_1  \Omega^{-1} (x(t))\]
according to the conjugacy relation.
By  \eqref{eq:Chazylike}  the  momentum component of $x(t)$ limits to $\cM v^+ (x(0))$ as $t \to \infty$.
On the other hand, by part 3  of the theorem we are  proving - the asymptotic near identity part, 
see \eqref{Moe:E:p},   the map $\Omega^{-1}$ tends to the identity along
escape orbits such as $x(t)$:   
$$\Omega^{-1} (x(t)) = x(t) + o(1) , \text{ as }  t \to \infty.$$  
Indeed, the term $q_{\min} ^{\alpha}$ appearing in estimate \eqref{Moe:E:p} tends to zero like
$t^{-\alpha}$ as $t \to \infty$. It follows that   
$\lim_{t \to \infty} {\rm pr}_1  \Omega^{-1} (x(t)) =  \cM v^+ (x(0))$,
which establishes \eqref{eq:meaningOmega}.

If $\Omega^{-1}$  mapped  {\it onto} $F_0$,  then the surjectivity of our map  \eqref{AFFINE}  would be immediate. 
$\Omega^{-1}$  would  map  $\cP_v$  onto the space of lines parallel to $v$ according to \eqref{eq:meaningOmega} and the conjugacy relation. 
And   $\Omega$, being the  inverse of $\Omega^{-1}$,
would be well-defined with domain all of $F_0$ and would map straight lines onto asymptotically free trajectories lying in $F^+$.
We could take $x^{(0)}(t)$ to be $\Omega(\ell_0 (t))$ where $\ell_0 (t))= (\cM v,  vt)$ corresponds to $b =0$.
Any $x(t) \in \cP_v$ can  be written, up to translation, uniquely  as $\Omega (\cM v,  vt + b)$ for some $b \in v^{\perp}$.  Moreover, both $\Omega$ and $\Omega^{-1}$ 
tend to the identity along escape orbits so that the limit in  \eqref{AFFINE} is the same as 
the limit achieved using the free flow, and so would yield $b = b(x, x^{(0)})$, and completing the proof.

$\Omega^{-1}$ is onto $F_0$  for non-singular potentials $V$.
To see this fact, observe that we can, in the  case of a non-singular potential, form   $\Phi_{-t} (x_0)$ for any $t$ and any $x_0$.
Incompleteness of the backward flow due to  collisions was the only thing  which prevented  
the direct Dollard-M\o ller map $\Omega$, defined  as the limit ${\Phi}_{-T} \circ \Phi^{T}_{T,0}$ as 
$T \to \infty$,  from existing and having domain all of $F_0$. 
The   analysis we    used  in  part 1 of the current theorem to insure  the existence of  $\Omega^{-1}$,  defined as  the limit of  $\Phi^{D}_{0,T}\circ\widehat{\Phi}_{T}$ as $T\to \infty$,
carries through  essentially verbatim  to yield the existence of 
 $\Omega: F_0 \to F^+$ and that it is   the inverse of our $\Omega^{-1}$.

We  deal with the case of singular potentials by observing that $\Omega^{-1}$
does not actually have to be onto, but only onto {\it modulo the  flow}, in order for the argument two paragraphs above to work.
For any $v_* \notin \Delta$ and    $b \in v_* ^{\perp} $   form $(v_*, b)$ and  denote its forward Dollard orbit by 
  $\Phi^D _{t,0} (v_*, b):= y_D (t;b), t_* \le t < \infty$.   
Eventually, for $t$ large enough,  we will show that these  Dollard orbits lies in the image of $\Omega^{-1}$, 
and that moreover $\Omega^{-1}$ is invertible there.  Then the entire Dollard ray $y_D ([t_*, \infty); b)$ will
lie in the image of $\Omega^{-1}$ and this will be enough.    
 To this end, fix any  relatively compact neighborhood $K$  of the origin in the full phase space $\bR^{dn} \times\bR^{dn}$.
Then there is a  $t$ large enough so that $K_t: =y_D (t) + K  \subset F^+_{\rm loc} \subset  F^+$.    To see this, 
observe that as $t$ increases without bound  the estimates
of \ref{def:finally:free}  must eventually   hold since    the $q_i$ occurring in the estimate   are equal to  $t v_*, i$ to leading order
while the $v_i$ are $v_{*, i}$.  It  thus follows from Theorem \ref{thm:mainMoller}  that $K_t  \subset F^+$
for  all sufficiently   large $t$.     Now,  as we saw a few paragraphs above,  part 3 (just proved; see also \eqref{Moe:E:p}) 
 tells us that the  map $\Omega^{-1}$ on $K_t$ is of the form $Id + h_t$ with $h_t =o(1)$ as $t \to \infty$. 
   As soon as $t$ is large enough so that the $C^{k-1}$-norm   of $h_t$ on $K_t$ is   less than $1$ we have that
$\Omega^{-1}$ is invertible  and that   $y_D (t) \in \Omega^{-1} (K_t) \cap K_t$.  We can let $t$ increase since the estimates only get better
and in this way conclude that  
the entire future Dollard ray $y_D ([t_*, \infty); b)$ lies in the image of  $\Omega^{-1}$, for some $t_* = t_* (v_*, b)$. Also  $\Omega$,
the inverse of $\Omega^{-1}$  exists along the Dollard ray.    
This analysis applies to any $b$, including $b = 0$. Now  take $v = v_*$
Take  for the `origin' of our trajectory space $\cP_v$ the solution $x^{(0)}(t) = \Omega (y_D (t; 0)$.
Then  
\[b(x, x^{(0)}) = \lim_{t \to \infty} q(t) - q_0 (t) = \lim_{t \to \infty} {\rm pr}_2 (y_D (t;b) - y_D (t; 0) )= b,\] where ${\rm pr}_2 (p, q)  = q$ is the projection onto configuration space.
We have proved that the map \eqref{AFFINE} is onto.    

Finally, to see that the map $x \mapsto b(x, x^{(0)})$ of \eqref{AFFINE} determines
the trajectory $x$ up to time translation use that fact that the same map determines the
trajectory up to time translation over on the free side, and that the free and Newtonian  limits are equal
since $\Omega^{-1}$ tends to the identity along escaping  orbits. 
\hfill $\Box$
 
\noindent

\begin{remark}
See also the remark in {\sc Derezi\'{n}ski} and {\sc G{\'e}rard} \cite[p.\ 24]{DG} regarding the affine structure of the tangent space and part 4.  
\end{remark} 

\begin{remark}[Derivation of asymptotics \eqref{eq:Chazylike}]
\label{rmk:asymptotics}

Set $\Phi^D _t  = \Phi^D_{t,0}$ so that $\Omega^{-1} = \lim_{t \to \infty} (\Phi^D _t)^{-1} \circ \Phi_t$.
It follows that for $t$ large we have $\Phi^D _t \circ  \Omega^{-1} = \Phi_t + o(1)$.
But $\Omega^{-1}$ tends to the identity along escaping orbits $x(t)$.  
This yields\\ 
$\Phi^D _t (x(T)) = \Phi_t (x(T)) + o(1)$
for $T$ sufficiently large,  $t \to \infty$ which is \eqref{eq:Chazylike}.   
\hfill $\Diamond$
\end{remark}  

\begin{remark}[Homogeneous potentials]\quad\\ 
If $V$ is a ($-\alpha$)--homogeneous potential, then
for every $k\in \bN$, $V$ is an $(\alpha,k)$--potential.
So in particular the Dollard-M\o ller transformation is smooth.
\hfill $\Diamond$
\end{remark}
\begin{earlier}\quad\\[-6mm] 
\begin{enumerate}[1.]
\item 
As we have indicated above, the assumption $\alpha\in (1/2,1]$ in Theorem \ref{thm:Dollard:Moeller}
can be relaxed, by generalizing the two-body technique from {\sc Herbst} \cite{He}. 
The price to be payed is a Dollard dynamics that is more involved than
\eqref{Dollard:dynamics}. 
\item 
Theorem 1 of {\sc Saari} \cite{Sa} states for the gravitational $n$--body system that under a
non-oscillation assumption the centers of mass of clusters asymptotically either move like 
$t\mapsto vt+D\log(t)+o(\log(t))$, or their mutual distances are of order $\cO(t^{2/3})$.
As this allows for non-trivial clusters, Saari's result is not contained in the statement of 
Theorem \ref{thm:Dollard:Moeller}.
On the other hand, Theorem \ref{thm:Dollard:Moeller} concerns general long range potentials
and controls the asymptotics of the flow, not just of individual orbits.
\item 
As Lemma \ref{lem:orbits:with equal:as:ve} below
shows, orbits with equal asymptotic momentum $\ov{p}^+$ synchronize their relative positions,
although their momenta $\tilde{p}$ approach $\ov{p}^+$ only slowly 
($\,\tilde{p}(t)-\ov{p}^+=\cO(t^{-\alpha})\,$).
See also {\sc Herbst} \cite[Lemma II.2]{He} for the case of potential scattering.
\hfill $\Diamond$
\end{enumerate}
\end{earlier}
\begin{lemma}[Orbits with equal asymptotic velocity]\label{rem:equal:as:vel}\quad\\ 
\label{lem:orbits:with equal:as:ve}%
For a long range potential $V$, consider initial conditions 
$x^{(i)}_0\equiv\big(p^{(i)}_0,q^{(i)}_0\big)\in\widehat{F}^\pm$ \ $(i=1,2)$,
whose asymptotic momenta $\ov{p}^+\big(x^{(i)}_0\big)$ respectively $\ov{p}^-\big(x^{(i)}_0\big)$ coincide.
Then
\beq
\Omega^{-1, \pm} \big( x^{(2)}_0 \big) - \Omega^{-1, \pm} \big( x^{(1)}_0 \big)
\ =\ \Big(0\ ,\ \lim_{t\to\pm\infty} \big(q(t,x^{(2)}_0)-q(t,x^{(1)}_0)\big) \Big).
\Leq{inverse:moeller:qminq} 
In particular, the limit on the right in~\eqref{inverse:moeller:qminq} is finite when the $x^{(i)}_0$ yield solutions having  the same asymptotic
velocities (or momenta).
\end{lemma}
\begin{remark}[Difference between long range and short range case]\quad\\ 
Note that the limits
$\lim_{t\to\pm\infty}\Big( \Phi_t \big(x^{(i)}_0 \big) - \Phi^D_{t,0}\circ \Omega^{*,\pm} \big(x^{(i)}_0 \big)\Big)$
do {\em not} exist for $(-\alpha)$--homogeneous potentials and $\alpha\in (0,1)$, see 
Appendix \ref{app:Dollard-Herbst}.
\hfill $\Diamond$
\end{remark}
\textbf{Proof of Lemma \ref{lem:orbits:with equal:as:ve}:}
With $a^\pm$ from \eqref{a:pm} and $\Omega_t$ from \eqref{Omega_T} we have
\begin{align}
(0,a^\pm)
&\stackrel{(1)}{=}
\lim_{t\to\pm\infty}  \big( p(t,x^{(2)}_0) - p(t,x^{(1)}_0) \, ,\, q(t,x^{(2)}_0) - q(t,x^{(1)}_0)\big)\NN\\
&\stackrel{(2)}{=}
\lim_{t\to\pm\infty} 
\big( \Phi^D_{t,0}\circ \Omega^{-1}_t(x^{(2)}_0)-\Phi^D_{t,0}\circ \Omega^{-1}_t(x^{(1)}_0) \big)\NN\\
&\stackrel{(3)}{=}
\lim_{t\to\pm\infty} \big(\Omega^{-1}_t(x^{(2)}_0) - \Omega^{-1}_t(x^{(1)}_0) \big)
\stackrel{\rm def.}{=}
\Omega^{-1, \pm}\big(x^{(2)}_0\big) - \Omega^{-1, \pm}\big(x^{(1)}_0\big)\, ,\NN
\end{align}
since\\
$\bullet$
By assumption $\ov{p}^\pm\big(x^{(i)}_0\big) = \lim_{t\to\pm\infty} p\big(t,x^{(i)}_0\big)$ coincide, and by
Lemma \ref{lem:q:minus:q} $a^\pm= \lim_{t\to\pm\infty} \big(q(t,x^{(2)}_0)-q(t,x^{(1)}_0)\big)$ 
exists, proving (1).\\
$\bullet$
Identity (2) follows from 
$\Phi^D_{t,0}\circ \Omega^{-1}_t=\Phi^D_{t,0}\circ\Phi^D_{0,t}\circ\widehat \Phi_t=\widehat \Phi_t$, see
\eqref{groupoid}.\\
$\bullet$
The Dollard dynamics  $\Phi^D$, see \eqref{Dollard:dynamics}, 
does not change momentum, which implies equality of the first components in (3).
Concerning the second components, 
\[Q\big(-t,p\big(t,x^{(i)}_0\big)\big) = q\big(t,x^{(i)}_0\big) 
- \cM^{-1}p\big(t,x^{(i)}_0\big) t - \int_0^t \nabla_p V\big(\langle s\rangle \cM^{-1}p\big(t,x^{(i)}_0\big)\big)\,ds\]
and $\big\| p\big(t,x^{(2)}_0\big) - p\big(t,x^{(1)}_0\big) \big\| = \cO(|t|^{-1-\alpha})$, see \eqref{pp:qq}. So
\[Q\big(\!-t,p\big(t,x^{(2)}_0\big)\big) - q\big(t,x^{(2)}_0\big) =
Q\big(\!-t,p\big(t,x^{(1)}_0\big)\big) - q\big(t,x^{(1)}_0\big) +\cO(|t|^{-\alpha}),\]
proving (3).
\hfill $\Box$\\[2mm]

Finally we prove a property special to ($-1$)--homogeneous potentials:
the existence of the Dollard-M\o ller transformation {\em and} of
asymptotes.  This property does not extend to $(1,k)$--potentials or
to ($-\alpha$)--homogeneous potentials, $0 < \alpha < 1$, as
counterexamples on the half-line show.

\begin{prop}[Asymptotes for ($-1$)--homogeneous potentials]\quad\label{prop:m:one:homog}\\
  Let $V$ be a $(-1)$--homogeneous potential and
  $\Phi^D_{\bullet,\bullet}$ its
  Dollard flow \eqref{Dollard:dynamics}.\\
  Then for all initial conditions $x_0\in \widehat{F}^\pm$ there exist
  unique $X_0^\pm\in F_0$ with
  \[\lim_{t\to \pm\infty}\big(\Phi^D_{t,0}(X_0^\pm) -
  \widehat{\Phi}_t(x_0)\big) = 0. \]
  In fact, $X_0^\pm=\Omega^{*,\pm}(x_0)$. 
\label{-1_Dollard}
\end{prop}
\textbf{Proof:} \\
We show the result for the limit $t\to+\infty$ and omit the superscript $\pm$ of $X_0^\pm$.
To make clear that $\alpha=1$ is the unique power with the described property, we first allow
for $\alpha$--homogeneous potentials with $\alpha\in (1/2,1]$.\\
We set $x(t)\equiv \big(p(t),q(t)\big) := \widehat{\Phi}_t(x_0)$ and 
$X(t)\equiv \big(P(t),Q(t)\big) :=\Phi^D_{t,0}(X_0)$ for the Dollard flow
with initial conditions $X_0:= \Omega^{*,+}(x_0)$ and show that the limit
$\lim_{t\to +\infty}\big(X(t) - x(t)\big)$ exists iff $\alpha=1$.\\
For all $t\in \bR$ we have $P(t) = \ov{p}^+(x_0) = \lim_{t\to +\infty} p(t)$.
So we must consider 
\begin{align}
F(t) &:= 
Q(t)-q(t) = Q_0 + \ov{v}^+(x_0)t 
+ \int_0^t \! \nabla_{\ov{p}^+} V\big(\langle s \rangle \cM^{-1}\ov{p}^+(x_0)\big)\,ds - q(t)\NN\\
&= Q_0 + \ov{v}^+(x_0)t + f_\alpha(t)\cM^{-1} \nabla V\big(\ov{v}^+(x_0)\big) - q(t) \,,\NN
\end{align}
(See Appendix  \ref{ex:Hom:Dollard}.) Its time derivative equals for $t>0$
\begin{align}
\dot{F}(t) &= 
\ov{v}^+(x_0) +  \langle t \rangle^{-\alpha} \cM^{-1} \nabla V\big(\ov{v}^+(x_0)\big) - \dot{q}(t)\NN\\
&=\cM^{-1} \Big[ \int_t^\infty\!\! \nabla V\big(q(s)\big)\,ds 
+ \langle t \rangle^{-\alpha} \nabla V\big(\ov{v}^+(x_0)\big)\Big] \NN\\
&=\cM^{-1} \Big[ \int_t^\infty\!\! \nabla V\big(\ov{v}^+(x_0) s + \cO\big(s^{1-\alpha}\log(s)\big)\big)\,ds 
+ \langle t \rangle^{-\alpha} \nabla V\big(\ov{v}^+(x_0)\big)\Big] \NN\\
&=\cM^{-1} \Big[ \int_t^\infty\!\! \nabla V\big(\ov{v}^+(x_0) s \big)\,ds 
+ \langle t \rangle^{-\alpha} \nabla V\big(\ov{v}^+(x_0)\big)\Big] + \cO\big( t^{-2\alpha} \log(t) \big)\NN\\
&=\big[  \langle t \rangle^{-\alpha} -\alpha^{-1} t^{-\alpha}\big] 
\cM^{-1} \nabla V\big(\ov{v}^+(x_0)\big) + \cO\big(t^{-2\alpha}\log(t)\big). \NN
\end{align}
We used $(-\alpha)$--homogeneity of $V$ in the second to last equation. To avoid 
a distinction of cases, we kept an $\cO(\log(t))$ term, that is unnecessary if $\alpha\in (1/2,1)$.\linebreak

So if $\alpha=1$, then the first term is of order $ \cO\big(t^{-3}\big)$, and only in this case
$\lim_{t\to +\infty}\big(Q(t) - q(t)\big)$ exists.
Subtracting this limit from $Q_0$, if non-zero, gives the unique initial conditions of the
Dollard flow that yield an asymptote.

However, for $\alpha=1$ we have $\lim_{t\to +\infty}\big(Q(t) - q(t)\big)=0$:
We just proved that the difference of the momentum $p(t)$ (which equals the momentum component of
$\Omega_t^{*,+}(x_0)$) and of  $P(t)$ is of order $\cO\big(t^{-2}\log(t)\big)$. So the difference of the
positions of the time $t$ Dollard flow with initial conditions $(P_0,Q_0):=X_0=\Omega^{*,+}(x_0)$
and $(P_t,Q_t):=X_t:=\Omega_t^{*,+}(x_0)$ is
\begin{align}
Q(t)-\big(Q_t+  t&\cM^{-1} p(t)+f_1(t)\nabla_{p(t)}V(\cM^{-1} p(t)\big) \label{Diff:Q:DM:asymp}\\
&=\; \big[Q_0-Q_t\big]\,+\,\big[\ov{v}^+(x_0)-\cM^{-1}p(t)\big] t \NN\\
&\quad \
+\sinh^{-1}(t)\Big[\nabla_{\ov{p}^+(x_0)}V(\cM^{-1} \ov{p}^+(x_0)\big) - \nabla_{p(t)}V(\cM^{-1} p(t)\big)\Big].
\NN
\end{align}
By definition of the Dollard-M\o ller transformation
$\lim_{t\to\infty}[Q_0-Q_t]=0$, whereas
$\big[\ov{v}^+(x_0)-\cM^{-1}p(t)\big] t=\cO\big(t^{-1}\log(t)\big)$,
and
\[\sinh^{-1}(t)\Big[\nabla_{\ov{p}^+(x_0)}V(\cM^{-1}
\ov{p}^+(x_0)\big) - \nabla_{p(t)}V(\cM^{-1} p(t)\big)\Big]\! 
=\cO\big(\log(t)\cdot t^{-2}\log(t)\big).\]
So the difference \eqref{Diff:Q:DM:asymp} has limit zero. \hfill $\Box$

%----------------------------------------------------
\section{On the scattering relation and map}
%----------------------------------------------------

When we replace limits $t \to \infty$ by $t \to - \infty$ we arrive at
the analogous objects for backward time, such as
$$v^{-} (x_0) = \lim_{t\to -\infty} q(t;x_0)/t,$$ 
in definition \ref{def:free}.  In this way we arrive at the backward
time analogue of being ``free'', which is to be in the set
$$F^-:=  \{x \in P : \text{the   solution   through }   x  \text{
  is backward free} \} ,$$
and the backward M\" oller transform
$$\Omega_-  := \lim_{t \to - \infty} \Phi_{-t} \circ\Phi^0_{t}: P
\dashedrightarrow P.$$ 
If $x_0 \in F^- \cap F^+$ then both $v^+ (x_0)$ and
$v^{-} (x_0)$ are defined, which leads us to the {\em scattering
  relation} $\sim_s$ on $\bR^{dn} \, \backslash\, \Delta$ under which
$v^- \sim_s v^+$ if and only if there exists an
$x_0 \in F^- \cap F^+$ such that $v^- (x_0) = v^-$ and
$v^+ (x_0) = v^+$.  Borrowing from quantum mechanics, the
``S-matrix'' or scattering map is defined by
$$S\, :=\,  \Omega_+ ^{-1}  \circ \Omega_- .$$
$S$ takes an ``initial condition'' $(p_-, C_-)$ at time $t = -\infty$
to an $x_0 \in F^- \cap F^+$ and then takes this $x_0$ to the
$(p_+, C_+)$ at $t = + \infty$ to which its solution corresponds.
Observe that $p_{\pm} = \cM v^{\pm} (x_0)$ so that the projection of
the {\em graph} of the scattering map onto its momentum components
$p_-, p_+$ yields the scattering relation (times the mass matrix
$\cM$).  We will leave this work to future researchers or future
times.
 
\begin{remark}[Manifold at infinity]\quad\\
  For an alternate construction of the scattering map which is valid
  for long range potentials and in particular for the Newtonian
  potential, see \cite{Duignan}. In this version $S$ is defined by
  adding a manifold at infinity and identifying the asymptotic
  velocities $v_-$, $v_+$ with equilibrium points at infinity. 
  \hfill $\Diamond$
\end{remark}

\appendix

\bigskip\bigskip
\noindent
{\bf \Large Appendices}

\vskip .1cm

In these appendices we go over some aspects of Dollard flows and
the induced transformations in a more leisurely fashion. 
Hamiltonians for homogeneous
potential are computed in  Appendix \ref{ex:Hom:Dollard}.  There we get two interesting surprises:
first, that a solution to the Dollard dynamics admits an asymptote if
and only if the initial velocity $a = \cM^{-1}p$ is a central
configuration in the sense of celestial mechanics.  The second
surprise is the appearance of hypergeometric functions.
In Appendix \ref{app:Dollard-Herbst} we show that for 
$0 < \alpha < 1$ one can actually
define two Dollard dynamics.  One  admits asymptotes but no M\o
ller transformation. The other, essentially the one we use,
admits no asymptotes but does yield a M\o ller transformation.  For
$\alpha = 1$ these two are equal, and this happy coincidence gives the
method more power here.  See Proposition \ref{-1_Dollard}.

\section{Precursors to Dollard-M\o ller}
\label{sec:Precursors}

We explore two
alternatives to the M\o ller transformation, by way of examples.
\begin{example}[Kepler-M\o ller transformation for the $n$--center
  problem]  
  After regularization, the motion of a single particle in the
  $n$--center potential\linebreak
  $V(q):=-\sum_{k=1}^n \frac{Z_k}{\|q-s_k\|}$ with $Z_k\in\bR$ and
  $s_k\in\bR^3$ leads to a complete, smooth flow $\Phi$.  By comparing
  the n-center flow with the regularized flow $\Phi^{(K)}$ of the
  Kepler Hamiltonian
  $H^{(K)}(p,q):=\eh \|p\|^2-\frac{Z_\infty}{\|q\|}$ with
  $Z_\infty:=\sum_{k=1}^n Z_k$, we can define a modified M\o ller
  transformation which exist for all initial values $x$ with
  $H^{(K)}(x)>0$, and which is smooth.  See \cite[section
  6]{Kn}.\hfill $\Diamond$%
\end{example}
In the case of the gravitational $n$--body problem we do not know of
any {\em time-independent} comparison Hamiltonian dynamics which
yields an explicit integrable flow and also yields a well-defined M\o
ller transformation.

\begin{example}[Asymptotes and Galilean boosts] \label{asymptotes}\quad\\%
  A positive {\em asymptote} for a solution curve $t \mapsto q(t,x)$
  to Newton's equations is, by definition, an affine line
  $L_+ \subseteq \bR^{dn}$ in configuration space whose distance
  $ \min_{\tilde{q}\in L_+} \|q(t)-\tilde{q}\|$ to the solution curve
  vanishes as $t\to +\infty$.  In a similar manner we define a
  negative asymptote $L_-$ by insisting its distance to the solution
  goes to zero as $t\to -\infty$.  Assuming that the asymptotic
  velocities $v^{\pm}(x)$ of a solution exist and are not zero, then
  $\lim_{t\to\pm\infty} \|q(t,x)\|=\infty$, and if an asymptote $L$
  exists, it is necessarily unique.  It is often the case that the
  asymptotes exist:
\begin{enumerate}[$\bullet$]
\item Solutions corresponding to short range potentials $V$ in the
  free region $F^+$ (definition~\ref{def:free}) always have
  asymptotes. This follows from the existence of inverse M\o ller
  transformation, proven in Theorem \ref{thm:both:moeller}.3.
\item 
Although the Kepler potential (or Newtonian 2-body problem)  is not   short range,
in the center of mass coordinates every Kepler hyperbola has  asymptotes in both directions.
\item
Similarly, a particle moving along a bi-hyperbolic orbit under the influence of the gravitational
or electrostatic potential due to  $n$ centers admits both positive and negative asymptotes,  
see \cite[section 6]{Kn}.
\end{enumerate}
The space of oriented affine lines in Euclidean $\bR^k$ is naturally
diffeomorphic to the cotangent bundle of the sphere $S^{k-1}$, and in
particular forms a symplectic manifold.  For an $n$-body problem in
$d$-dimensions, we have $k =nd$ and may try to construct a substitute
for the M\o ller transformation by sending
$L_- (x_0) \mapsto L_+ (x_0) \in T^* S^{k-1}$,  where $L_{\pm} (x_0)$
are the positive and negative asymptotes of the initial condition
$x_0$.  With some luck, $L_{\pm} (x_0)$ might exist for all
$x_0 \in F^+ \cap F^-$ and we would then have our scattering map
as a map between open subsets $T^*S^{k-1}$, a symplectomorphism even
if the asymptotics of solutions depended sufficiently smoothly on the
initial condition $x_0$.  See for example \cite[section 16]{Kn}.

Galilean boosts destroy the existence of asymptotes, so that we cannot
expect asymptotes to exist for general long range potentials.  To see
this destruction phenomenon, take for simplicity $n =2$,  $d =1$ and
$m_1 = m_2 = 1$.  The Newtonian two-body equations read
$\ddot q_1 = (q_2 - q_1 )/ r^3$, $\ddot q_2 = (q_2 -q_1)/ r^3$ with
$r = |q_1 -q_2|$. \\ 
The reduced mass is $ \mu = 1/2$ so that the
corresponding Kepler problem becomes
$ \ddot x = - 2x/|x|^3, x = q_1 - q_2$.  A hyperbolic solution
$q(t) = (q_1 (t), q_2 (t))$ with asymptotic energy $1$ in the center
of mass frame $q_1 + q_2 = 0$ will have asymptotics
$q_1 (t) = t + 2 \log(t) + c + o(1)$, $q_2 (t) = -t - 2 \log(t) - c + o(1)$.\\
Apply the Galilean boost $(q_i, t) \mapsto (q_i + t v, t)$ to this
solution to obtain a new two-body solution
$\tilde q = (\tilde q_1, \tilde q_2)$ whose asymptotic expansion is
\begin{align*}
\tilde q_1 (t) &=\ \ t + vt  +2  \log(t) + c + o(1),\\
\tilde q_2 (t) &= - t + vt   -2  \log(t) - c + o(1).
\end{align*}
Now, the signed distance between a  point $Q = (x,y) $  and a line $L \subseteq \bR^2$ 
is given by  the affine expression 
 $d(Q, L) = u \cdot Q + e = ax + by + e$ where $u=(a,b)$ is the  unit vector perpendicular to the
direction of $L$ and where  $e$ is the distance between $L$ and $(0,0)$.   Consequently
the signed distance between  our putative $L$ and our moving solution 
$(\tilde q_1 (t), \tilde q_2 (t))$ must have the form
\[a( t+ v t +2 \log (t)) + b (-t + vt  -2 \log(t))  + e  + o(1).\]
Expanding out we find that this signed distance has asymptotic expansion\\  
$(a - b + 2v ) t + 2 (a- b)  \log(t)  + e +  o(1)$.
For this distance to tend to zero with$t$   we must have that $a-b +2v  = 0$
as well as  $a-b = 0$ which is impossible if $v \ne 0$.~$\Diamond$%
\end{example}
The latter example shows that we cannot use asymptotic affine lines to
model scattering for long range potentials. We do not know how to use
a time independent model flow as a replacement to free flow, or the
Kepler flow of the previous example either.   

%----------------------------------------------------
\section{Dollard and central configurations for homogeneous potentials} 
\label{ex:Hom:Dollard}
%----------------------------------------------------
 
Take $V$ to be one of the power law potentials of homogeneity
$-\alpha$, $\alpha\in(1/2,1]$, as defined by \eqref{hom:pot}.  These
potentials are strictly long range, that is, they are long range and
not short range.\footnote{The $(-\alpha)$--homogeneous potentials
  with $\alpha>1$ are of short range. So the time-inde\-pen\-dent
  kinetic Hamiltonian $K$ can be used for defining the M\o ller
  transformations.}  The main interest is of course in the case
$\alpha=1$, including gravitational or electrostatic interactions.

We first note that, although the potential $V$ is unbounded in our
case, the statements of Theorem \ref{thm:both:moeller} concerning
long-range pair interactions apply, since for any $k\in\bN$ the norms
\eqref{V:k:norm} are still finite for these homogeneous potentials.

The proper Dollard Hamiltonian $V \big( \langle t\rangle\cM^{-1}p \big)$ equals
\beq
\tilde{H}^D_t(p,q) =\LA t\RA^{-\alpha}\sum_{1\le i< j\le n} \frac{ I_{i,j}}{\|v_i - v_j\|^\alpha} \qquad(t\in \bR),
\Leq{dollard:hamiltonian:homogeneous}
and for $W$, defined in \eqref{def:W}, we explicitly get 
\begin{align}
\int_0^t \nabla_p V \big(\langle s\rangle \cM^{-1}p \big)\,ds
&= \int_0^t \langle s\rangle^{-\alpha} \,ds \ \nabla_p V(\cM^{-1}p)
= f_\alpha(t)\ \nabla_p V(\cM^{-1}p),\NN
\end{align}
with $f_\alpha\in C^\infty(\bR,\bR)$, 
$f_\alpha(t):=t \; _2F_1\left({\textstyle \frac{1}{2},\frac{\alpha}{2};\frac{3}{2};-t^2}\right)$ being odd,
and $_2F_1$ denoting the hypergeometric function.
For $\alpha=1$ this simplifies to $f_1(t)=\sinh ^{-1}(t)$.
Moreover,
\[  \nabla_p V(\cM^{-1}p) = \textstyle \Big(\sum\limits_{k\in N\setminus\{1\}} 
\frac{-\alpha I_{k,1}(v_k-v_1)}{m_1\|v_k - v_1\|^{\alpha+2}}
\; ,\ldots,\; \sum\limits_{k\in N\setminus\{n\}} 
\frac{-\alpha I_{k,n}(v_k-v_n)}{m_n\|v_k - v_n\|^{\alpha+2}}\Big) . \]

\begin{definition}
A vector $x\in \bR^{dn} \,\backslash\,\Delta$ is called a {\em central configuration} if it is 
linearly dependent with respect to $\cM^{-1}\nabla V(x)$.
\end{definition}
We see from the definition \eqref{def:W} of $W$ and the last formulae above 
that for these $(-\alpha)$--homogeneous potentials the Dollard dynamics has an asymptote
in the sense of Example  \ref{asymptotes} 
if and only if $v$ is a central configuration.\footnote{This has been noted independently
by Alain Albouy (private communication with A.K.)}
%

%----------------------------------------------------
\section{Dollard and the two methods of Herbst} \label{app:Dollard-Herbst}
%----------------------------------------------------

The point of this appendix is to show that a Dollard dynamics cannot 
both lead to asymptotics for the scattering solution {\em and}
existence of M\o ller transformations. Our  work here  follows the
ideas of  Herbst \cite{He}.

For $\alpha=1$ the Hamiltonian \eqref{dollard:hamiltonian:homogeneous}
is a classical analogue of the quantum ansatz introduced by {\sc
  Dollard} in \cite{Do}.  For the potential scattering of two
particles, {\sc Herbst} translated Dollard's ideas to the classical
case and generalized them to long range potentials in his interesting
article \cite{He}.  See also \cite[section 1.12]{DG}.

Herbst actually describes {\em two} different natural definitions of a
Dollard Hamiltonian for $(-\alpha)$--homogeneous potentials, when
$\alpha\in(0,1)$.  We proceed to describe and analyze these in the
one-dimensional case  ($d =1$).

So   consider   
$H:T^*\bR^+\to\bR$, $H(p,q):=\eh p^2+V(q)$ with $V(q):=I/q^\alpha$ the Hamiltonian flow.  
For initial conditions $x_0=(p_0,q_0)$ with velocity $p_0>0$ and $h:=H(p_0,q_0)>0$
the asymptotic velocity $p^+$ equals $\sqrt{2h}$. By assuming $\alpha\in(1/2,1)$, 
we avoid the necessity of multiple iterations of integral equations, which would 
only blur the basic phenomenon.

As $\ddot{q}=\alpha I q^{-1-\alpha}$, in the large time asymptotics the solution has the form
\beq
q(t;p_0,q_0) = p^+t -\frac{I  \l((q_0+p^+t)^{1-\alpha}-q_0^{1-\alpha}\ri)}{(1-\alpha)(p^+)^2}
 + q_0 + \delta q_\alpha(p^+,q_0) + \cO\big(t^{1-2\alpha}\big),
\Leq{true:sol}
for $\alpha\in(1/2,1)$, respectively
\[q(t;p_0,q_0) = p^+t - \frac{I}{(p^+)^2}\log(1+p^+t/q_0) + q_0 + \delta q_1(p^+,q_0) + \cO(t^{-1}\log(t))\]
for $\alpha=1$, 
with $\lim_{q_0\to\infty}\delta q_\alpha(p^+,q_0)=0$.
\begin{enumerate}[$\bullet$]
\item 
We apply {\em Herbst's second method}, leading to a Dollard type
M\o ller transformation (Theorem  III.1 of \cite{He}). Thus we obtain a sequence
of time dependent Hamiltonians $H_t^{(k)}(p):=\eh p^2 +U^{(k)}(p,t) $ for
\[ U^{(0)}(p,t):=0\qmbox{,}U^{(k+1)}(p,t) := V\l(pt+\textstyle\int_0^t D_1U^{(k)}(p,s)\,ds\ri),\]
independent of $q$ but dependent on the asymptotic velocity $p$.
So
\[ U^{(1)}(p,t)=\frac{I}{(pt)^\alpha}
\qmbox{,} U^{(2)}(p,t) =  
I \left(p t \left(1-\frac{\alpha I (p t)^{-\alpha}}{(1-\alpha) p^2}\right)\right)^{-\alpha}.\]
For $\alpha\in (1/2,1)$ the solutions of the Hamiltonian equations for $H^{(k)}$ are 
\beq
q^{(0)}(t;p^+,q_0) =  p^+t + q_0  \qmbox{,}
q^{(1)}(t;p^+,q_0) = p^+t -\frac{\alpha}{1-\alpha} \frac{I}{(p^+)^{1+\alpha}} t^{1-\alpha}  + q_0\, ;
\Leq{Herbst:3}
$q^{(k)}$ for $k\ge2$ give corrections to $q^{(1)}$ with negative asymptotic order in $t$.

Now if one compares $q^{(1)}$ in \eqref{Herbst:3} with the asymptotics 
\eqref{true:sol} of the true solution, then one notices in
the term asymptotic to a multiple of $t^{1-\alpha}$ an additional factor~$\alpha$.  

Without much calculation, one sees that the inverse M\o ller transform 
$(\Omega^+)^{-1}$ exists
for initial condition $x_0:=(p_0,q_0)$ with $q_0:= 0_+$:
Then by  \eqref{true:sol} the solution of the initial value problem with Hamiltonian $H$ equals
\begin{align}
p(t; x_0)&= 
p^+-\frac{I}{(p^+)^{1+\alpha}}t^{-\alpha}+\cO(t^{-2\alpha})\,,\NN\\
q(t; x_0)&=p^+t -\frac{I}{(1-\alpha)(p^+)^{1+\alpha}}t^{1-\alpha}+\delta q+\cO(t^{1-2\alpha})\, .\NN
\end{align}
Note that by our assumption $\alpha\in(1/2,1)$ the function $t\mapsto t^{-2\alpha}$ is in $L^1\big([1,\infty)\big)$.
For $x_1:=(p_1,q_1)$ we have by \eqref{Herbst:3}
\beq
q^{(1)}(-t;x_1) = - p_1 t + \frac{\alpha I} {(1-\alpha)p_1^{1+\alpha}}t^{1-\alpha}+q_1.
\Leq{Moeller:convergent}
Setting $x_1:=\big(p(t; x_0),\, q(t; x_0)\big)$, we get convergence
of \eqref{Moeller:convergent} as $t\to+\infty$.
\item
We now apply {\em Herbst's first method}, leading to a solution
whose difference to $q(t;p_0,q_0)$ converges as $t\to\infty$ (Theorem  II.1 of \cite{He}). 
So we iteratively define functions $z^{(k)}$
of time $t$ and asymptotic momentum $p$ by setting
\[z^{(k)}(0,p):=0 \mbox{ , }\ \dot{z}^{(0)}(t,p) := p \qmbox{and}
\dot{z}^{(k+1)}(t,p) := p - \int_t^\infty\!\! F\big(z^{(k)}(s,p)\big)\,ds,\]
with force $F(q):=-\nabla V(q)= \alpha Iq^{-1-\alpha}$.
We obtain $z^{(0)}(t,p^+)=p^+t$ and 
\beq
z^{(1)}(t,p^+) = p^+t - \frac{1}{1-\alpha} \frac{I}{(p^+)^{1+\alpha}} t^{1-\alpha}.
\Leq{Herbst:2}
If we compare $z^{(1)}$ from \eqref{Herbst:2} with \eqref{true:sol}, we see that 
the factors of the terms diverging as $t\to+\infty$ agree. So here 
the solution \eqref{true:sol} has a limit
$\lim_{t\to\infty} \big(q(t;p_0,q_0)-z^{(1)}(t,p^+) \big)$; 
the same is true for the time derivatives.
However, for \eqref{Herbst:3} the corresponding limit does {\em not} exist if $\alpha<1$.
\end{enumerate}
Due to the appearance of the regularization $\langle t\rangle$ of $|t|$,
the Hamiltonian dynamics generated by our Dollard Hamiltonian \eqref{dollard:hamiltonian:homogeneous}
does not coincide with the ones of Herbst's first or second method. \hfill $\Diamond$

So we have seen that for general long range potentials a Dollard
dynamics cannot both lead to asymptotics for the scattering solution
{\em and} M\o ller transformations.

A comparison of the $\alpha$--dependence of 
$q^{(1)}$ in \eqref{Herbst:3} and of $z^{(1)}$ in \eqref{Herbst:2} as $\alpha\!\nearrow \!1$ suggests that 
both properties could coincide for $\alpha=1$, the Kepler potential. This is indeed the case, 
as we see in the body of the paper,  in Proposition \ref{prop:m:one:homog}.

%-------------------------------------------------------------------------------------------------------
%
\addcontentsline{toc}{section}{References}
%

%-------------------------------------------------------------------------------------------------------
\end{document}